\documentclass[onefignum,onetabnum]{siamart190516}
\usepackage{animate}
\usepackage{amsmath}
\usepackage{amssymb}
\usepackage{graphics}
\usepackage{graphicx}
\usepackage{footnote}
\usepackage{textcomp}
\usepackage{mathrsfs}
\usepackage{epstopdf}
\usepackage{array}
\usepackage[maxfloats=99]{morefloats}
\usepackage{url}
\usepackage{cases}
\usepackage{mathscinet}
\usepackage[normalem]{ulem}
\usepackage[noend]{algpseudocode}
\usepackage{color}
\usepackage{cite}
\usepackage{bm}
\usepackage{subcaption}

\usepackage{stmaryrd}
\usepackage{color}
\usepackage{mathtools}
\usepackage{cleveref}

\graphicspath{{./fig/}{./fig/dr/}{./fig/other/}}

\makeatletter
\def\BState{\State\hskip-\ALG@thistlm}
\makeatother

\numberwithin{equation}{section}

\newtheorem{remark}{Remark}[section]

\usepackage{xparse}

\NewDocumentCommand{\dgal}{sO{}m}{%
  \IfBooleanTF{#1}
    {\dgalext{#3}}
    {\dgalx[#2]{#3}}%
}

\NewDocumentCommand{\dgalext}{m}{%
  \sbox0{%
    \mathsurround=0pt 
    $\left\{\vphantom{#1}\right.\kern-\nulldelimiterspace$%
  }%
  \sbox2{\{}%
  \ifdim\ht0=\ht2
    \{\kern-.45\wd2 \{#1\}\kern-.45\wd2 \}%
  \else
  \fi
}

\NewDocumentCommand{\dgalx}{om}{%
  \sbox0{\mathsurround=0pt$#1\{$}%
  \sbox2{\{}%
  \ifdim\ht0=\ht2
    \{\kern-.45\wd2 \{#2\}\kern-.45\wd2 \}%
  \else
    \mathopen{#1\{\kern-.5\wd0 #1\{}
    #2
    \mathclose{#1\}\kern-.5\wd0 #1\}}
  \fi
}

\DeclareFontEncoding{FMS}{}{}
\DeclareFontSubstitution{FMS}{futm}{m}{n}
\DeclareFontEncoding{FMX}{}{}
\DeclareFontSubstitution{FMX}{futm}{m}{n}
\DeclareSymbolFont{fouriersymbols}{FMS}{futm}{m}{n}
\DeclareSymbolFont{fourierlargesymbols}{FMX}{futm}{m}{n}
\DeclareMathDelimiter{\VERT}{\mathord}{fouriersymbols}{152}{fourierlargesymbols}{147}

\usepackage[T1]{fontenc}
\usepackage{textcomp}
\usepackage[scaled]{helvet}

\title{Deep unfitted Nitsche method for elliptic  interface problems }
\author{Hailong Guo\thanks{School of Mathematics and Statistics,  The University of Melbourne,  Parkville, VIC 3010, Australia   (hailong.guo@unimelb.edu.au).}
\and %
Xu Yang\thanks{Department of Mathematics, University of California, Santa Barbara, CA, 93106, USA (xuyang@math.ucsb.edu).}
}

\begin{document}

\maketitle

\begin{abstract}
This paper proposes a deep unfitted Nitsche method for computing elliptic interface problems with high contrasts in high dimensions. To capture discontinuities of the solution caused by interfaces, we reformulate the problem as an energy minimization involving two weakly coupled components. This enables us to train two deep neural networks to represent two components of the solution in high-dimensional. The curse of dimensionality is alleviated by using the Monte-Carlo method to discretize the unfitted Nitsche energy function. We present several numerical examples to show the performance of the proposed method.


\end{abstract}



\section{Introduction}

In this paper, we continue our previous studies on elliptic interface problems \cite{guo2017gradient,guo2018gradient,GuYa2018}, arising in many applications such as fluid dynamics and materials science, where the background consists of rather different materials on the subdomains separated by smooth curves (or surfaces) called interfaces. We aim to address the high-dimensional challenge, which is well known as the curse of dimensionality leading to unaffordable computational time in traditional numerical methods ({\it e.g.}, finite difference and finite element methods).

Deep neural networks have been shown as a powerful tool to overcome the curse of dimensionality \cite{dudley1969speed,fournier2015rate,weed2019sharp,dos2018simulation}, and have been applied to solve partial differential equations (PDEs), {\it e.g.}, the deep BSDE method \cite{HJE2018, EHJ2017},  the deep Galerkin method (DGM)\cite{SiSp2018},  the physics-informed neural networks (PINNs) \cite{raissi2019physics},   the deep Ritz method (DRM) \cite{EYu2018}, and the weak adversarial networks (WAN) \cite{ZBYZ2020}. The deep BSDE reformulates the time-dependent equations into stochastic optimization problems. DGM and PINNs train neural networks by minimizing the mean squared error loss of the equation residual, while DRM trains networks by minimizing the energy functional of the variational problem equivalent to the PDEs. 
WAN uses the weak formulation and trains the primary and adversarial network alternatively using the min-max weak formulation.   Moreover,  the convergence of DRM was studied by \cite{lu2021priori,duan2021convergence}, and the deep Nitsche method was proposed in \cite{LiMi2021}, which enhanced the deep Ritz method with natural treatment of essential boundary conditions. In a recent work\cite{ShYa2021},  Sheng and Yang trained an additional neural network to impose the Dirichlet boundary conditions. 
However, in general, these neural network-based methods require the smoothness of the solutions to the PDEs. They thus can not be directly used to solve the elliptic interface problems, where the solutions are only piecewise smooth.

In literature, there are some recent works of solving elliptic interface problems using neural networks. For example, \cite{WaZh2020} proposed a network architecture similar to the deep Ritz method \cite{EYu2018}, and solved the equivalent variational problem with the boundary conditions approximated by a shallow neural network. \cite{he2020mesh} used different neural networks to approximate the solutions in disjoint subdomains. They reformulated the interface problem as a least-squares problem and solved it by stochastic gradient descent. \cite{hu2021discontinuity} proposed the discontinuity capturing shallow neural network (DCSNN) to approximate piecewise continuous functions and solved elliptic interface problems by minimizing the mean squared error loss consistent with the residual of the equation, boundary and interface jump conditions.

In this paper, we propose a deep learning method for interface problems based on the minimization of the unfitted Nitsche's energy, inspired by our previous studies \cite{GuYa2018,GYZ2021, GuYZ2021} on the unfitted Nitsche's method. One of the most significant differences between the unfitted Nitsche's method\cite{BCHLM2015, HaHa2002,GYZ2021, GuYZ2021} and other methods for interface problems ({\it e.g.}, the immersed finite element method \cite{LiIt2006}) is that the unfitted Nitsche's finite element functions can be discontinuous inside elements. This is possible by adopting two different sets of basis functions on the interface elements ({\it i.e.}, the elements cut by the interface) which are weakly coupled together using Nitsche's methods.  Based on the unfitted Nitsche's formulation, we can define a so-call unfitted Nitsche energy functional (see equation \eqref{equ:energy} ). It turns out that the weak formulation of unfitted Nitsche's method is just the Euler-Lagrange equation of unfitted Nitsche energy functional.  To address the challenges of the curse of dimensionality, we naturally use deep neural networks to present functions in high dimensions. Following the idea of classical unfitted Nitsche's method \cite{BCHLM2015, HaHa2002, GYZ2021, GuYZ2021}, we use two deep neural networks: one for the part inside the interface and the other one for the region outside the interface. These two parts are weakly connected using Nitsche's method. The deep unfitted Nitsche method trains the two neural network functions independently using the same unfitted Nitsche energy functional.


The rest of the paper is organized as follows. In Section~\ref{sec:model}, we introduce the model setup of the elliptic interface problems and its unfitted Nitsche's weak form; In Section~\ref{sec:deepNitsche}, we describe the formulation of deep unfitted Nitsche method with details; We present numerical examples in Section~\ref{sec:example}, and make conclusive remarks in Section~\ref{sec:conclusion}.

\section{Model equation}\label{sec:model} Denote by $\Omega$ a bounded Lipschitz domain in $\mathbb{R}^d$, and assume there is a smooth closed curve (or surface)  $\Gamma$ separating $\Omega$ into two parts: $\Omega_1$ and $\Omega_2$. In general, $\Gamma$ can be described as a zero level set of some level set function $\phi$. Then we have $\Omega_1 = \{x\in\Omega |\phi(x)<0\}$ and $\Omega_1 = \{x\in\Omega |\phi(x)>0\}$.

\begin{figure}[!h]
   \centering
  \includegraphics[width=0.5\textwidth]{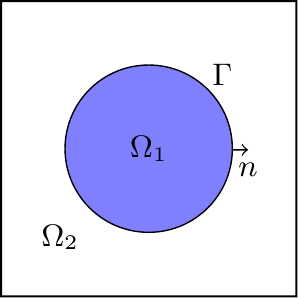}
   \caption{An illustrative example of a domain $\Omega$ with a circular interface $\Gamma$ in two dimension. Here $n$ stands for the outer normal direction of the inner domain $\Omega_1$. $\Omega_2$ is the outer domain.}
   \label{fig:circle}
\end{figure}

In this paper, we shall consider the following elliptic interface problem
\begin{subequations}\label{equ:interface}
 \begin{align}
  -\nabla \cdot (\beta(x) \nabla u(x)) &= f(x),  \quad \text{ in } \Omega_1\cup\Omega_2, \label{equ:model}\\
   u & = g, \quad\quad\,\,  \text{ on } \partial\Omega, \label{equ:bnd}\\
      \llbracket u\rrbracket &=p,  \quad\quad\,\,  \text{ on } \Gamma \label{equ:valuejump}\\
        \llbracket \beta \partial_n u \rrbracket&= q,   \quad\quad\,\,  \text{ on } \Gamma\label{equ:fluxjump}
\end{align}
\end{subequations}
where  $\partial_n u=(\nabla u)\cdot n$ with $n$ being the unit outward normal vector of $\Gamma$   and  the
jump $\llbracket w \rrbracket$ on $\Gamma$  is defined as
\begin{equation}\label{equ:jump}
\llbracket w\rrbracket = w_2|_{\Gamma} - w_1|_{\Gamma},
\end{equation}
with $w_i = w|_{\Omega_i}$ being the restriction of $w$ on $\Omega_i$.
The diffusion coefficient $\beta(x) $ is a piecewise constant, {\it i.e.},
\begin{equation}
\beta(x) =
\left\{
\begin{array}{ccc}
    \beta_1 &  \text{if } x\in \Omega_1, \\
    \beta_2  &   \text{if } x\in \Omega_2,\\
\end{array}
\right.
\end{equation}
which has a finite jump of function value at the interface $\Gamma$. 
Without loss of generality, we assume $\beta_0=\min(\beta_1,\beta_2) > 0$. An illustration of the domain $\Omega$ is given in Figure \ref{fig:circle}.

To prepare the presentation of Nitsche's weak formulation, we  introduce two  weights 
\begin{equation}
\kappa_1 = \frac{\beta_2}{\beta_1+\beta_2}, \quad \kappa_2 = \frac{\beta_1}{\beta_1+\beta_2},
\end{equation}
which satisfies that $ \kappa_1 + \kappa_2=1$.
Then,  we define the weighted  averaging of a function $w$ on the interface $\Gamma$ as
\begin{equation}\label{equ:avg}
\dgal{w} = \kappa_1w_1|_{\Gamma} +\kappa_2w_2|_{\Gamma}, \end{equation}
and also its dual weighted averaging as
\begin{equation}\label{equ:dualavg}
 \dgal{w}^{\ast}= \kappa_2w_1|_{\Gamma} +\kappa_1w_2|_{\Gamma}.
\end{equation}

Let $H^{1}(\Omega_1\cup \Omega_2)$  be the function space consisting of piecewise Sobolev  functions $w$  such
 that $w|_{\Omega_1}\in H^1(\Omega_1)$ and $w|_{\Omega_2}\in H^1(\Omega_2)$, whose  norm  is defined as
 \begin{equation}\label{eq:pnorm}
\|w\|_{1, \Omega_1\cup \Omega_2} = \left( \|w\|_{1, \Omega_1}^2 + \|w\|_{1, \Omega_2}^2\right)^{1/2},
\end{equation}
where $\|\cdot \|_{1, \Omega_i}$ is the $H^1$-norm of a function in $H^1(\Omega_i)$. Similar notations are used for piecewise $L^2$ space and its corresponding norm. 
In addition, let $H^{1}_g(\Omega_1\cup \Omega_2)$  be the subspace of $H^{1}(\Omega_1\cup \Omega_2)$ such that  $u|_{\partial\Omega_2}=g$. In particular, $H^{1}_0(\Omega\cup\Omega)$ is the subspace  of $H^{1}(\Omega_1\cup \Omega_2)$ with homogeneous Dirichlet boundary conditions. 

The unfitted  Nitsche's weak formulation \cite{HaHa2002, BCHLM2015, GuYa2018} of the interface problem \eqref{equ:model}-\eqref{equ:fluxjump} is to 
find $u \in H^{1}(\Omega_1\cup \Omega_2)$ such that
\begin{equation}\label{equ:var}
a(u, v) = \ell(v), \quad \forall  v \in H^{1}_0(\Omega_1\cup \Omega_2);
\end{equation}
where the bilinear form $a(\cdot, \cdot) $ is defined as
\begin{equation}\label{equ:bilinear}
\begin{split}
 a(u,v) = &\sum\limits_{i=1}^2\left(\beta\nabla u_{i}, \nabla v_{i}\right)_{\Omega_i}
- \left\langle  \llbracket u \rrbracket ,\dgal{\beta \partial_nv}\right\rangle_{\Gamma} \\
&-\left\langle  \llbracket v \rrbracket ,\dgal{\beta \partial_nu}\right\rangle_{\Gamma}
+ \gamma_f\left\langle\llbracket u \rrbracket,  \llbracket v \rrbracket\right\rangle_{\Gamma}, 
\end{split}
\end{equation}
and the linear functional $\ell(\cdot)$ is defined as
\begin{equation}\label{equ:functional}
\begin{split}
\ell(v) = &\sum\limits_{i=1}^2(f, v_{i})_{\Omega_i} - \langle p,\dgal{\beta\partial_n v}\rangle_{\Gamma}
+ \gamma_f\left\langle  p,  \llbracket v \rrbracket  \right\rangle_{\Gamma}
+\left\langle q, \dgal{v}^{\ast}\right\rangle_{\Gamma},
\end{split}
\end{equation}
with the stability parameter $\gamma_f$ and $\langle \cdot, \cdot \rangle_{\Gamma}$ being the $L^2$-inner product on $\Gamma$.

To adopt the deep neural network, we reformulate the variational problem \eqref{equ:var} as an energy minimization problem.  For such a purpose, we define the  unfitted Nitsche's  energy functional $L_n$ as 
\begin{equation}\label{equ:energy}
	\begin{split}
		L_n(v):= &\frac{1}{2}a(v,v) - \ell(v)\\
		= &\frac{1}{2}\sum_{i=1}^{2}\int_{\Omega_i}\beta_i\nabla v_i\cdot\nabla v_idx + \frac{\gamma_f}{2}\int_{\Gamma}\left(  \llbracket v \rrbracket -
		 p \right)^2ds + \\
		& \int_{\Gamma}(p-  \llbracket v \rrbracket)\dgal{\beta \partial_nv}ds
		- \int_{\Omega}fvdx - \\
		& \int_{\Gamma} q\dgal{v}^{\ast}ds -\frac{\gamma_f}{2}\int_{\Gamma} p^2ds.
	\end{split}
\end{equation}

Then, we can show the variational problem is equivalent to the following energy minimization problem:
 \begin{equation}\label{equ:energymin}
 	 u = \arg\min_{v\in H^1_g(\Omega)} L_n(v).
 \end{equation}
 In other words, equation \eqref{equ:model}-\eqref{equ:fluxjump} is the Euler-Lagrangian equation of \eqref{equ:energymin}.

\begin{remark}
To simplify the presentation of this paper, we only consider inhomogeneous Dirichlet boundary conditions. For Neumann and Robin boundary conditions, they can be absorbed into the variation formulation. Other types of boundary conditions like mixed Dirichlet and Neumann boundary conditions can be handled similarly. 
\end{remark}

 \begin{remark}
 	We choose to  impose the Dirichlet boundary condition strongly by building it into the solution space.  Alternatively, we can impose the Dirichlet boundary condition weakly as in \cite{LiMi2021 }. 
 \end{remark}
 
 \section{Deep unfitted Nitsche method}\label{sec:deepNitsche}
 In this section, we use the universal approximation property of deep neural network and flexibility of unfitted Nitsche weak form to develop the deep unfitted Nitsche method for solving elliptic interface problem \eqref{equ:model}-\eqref{equ:fluxjump}. In the first subsection, we briefly introduce the deep neural networks used in this paper. In the second subsection, we describe the details of the deep unfitted Nitsche formulation. 
 
\subsection{Deep neural network} One of the critical ingredients for using deep learning to solve partial differential equations is to select a deep neural network as the ansatz function for trial functions. The commonly used deep neural networks to approximate the solutions to PDEs include feedforward neural network\cite{LMMK, GoBC2016}, and residual neural network (ResNet)\cite{EYu2018, HZRS2016}. Similar to the deep Ritz method in \cite{EYu2018},  we choose the ansatz function to be the ResNet.

For each integer  $i$ and some positive integer $m$, let $W^{[i_j]}\in \mathbb{R}^{m\times m}$ be the matrix of weights and $b^{[i_j]}\in\mathbb{R}^m$ be vector of biases for $j=1,2$. Then, the $i$th block of ResNet with $m$ neurons can be written as 
\begin{equation}
f_i(s) = \sigma (W^{[i_2]} \sigma(W^{[i_1]}s + b^{[i_1]} ) +  b^{[i_2]}) + s,
\end{equation}
where $\sigma$ is the activation function \cite{GoBC2016, HiHi2019}. To avoid the problem of vanish gradient, smooth functions like sigmoid and hyperbolic tangent can be adopted.

\begin{figure}[!h]
   \centering
  \includegraphics[width=\textwidth]{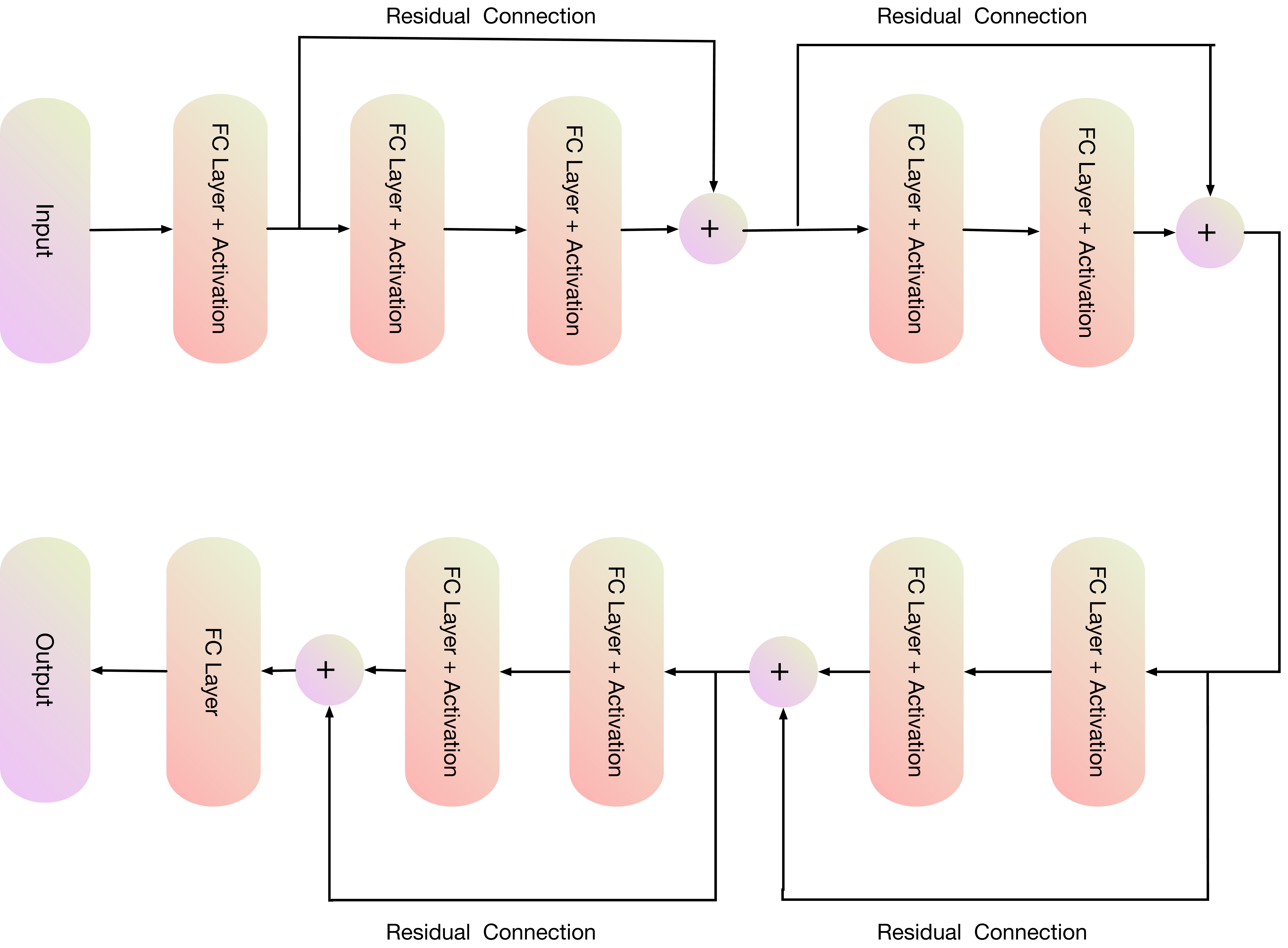}
   \caption{We present a diagram of ResNet with four blocks where FC layer denotes fully connected layer. Similar to the deep Ritz method in \cite{EYu2018},  we choose the ansatz function to be the ResNet.}
   \label{fig:resnet}
\end{figure}

Similarly, let $W^{[0]} \in \mathbb{R}^{m,d}$ (or $W^{[n+1]} \in \mathbb{R}^{1, m}$) be the matrix of weights in the first (or last) layer and $b^{[0]}\in\mathbb{R}^m$ (or $ b^{[n+1]}\in\mathbb{R}^1)$ be  the vector of biases in the first (or last) layer.
The ResNet with $n$ block can viewed as the composition of $f_i$'s  
\begin{equation}\label{equ:dnn}
	u_{\theta}(x) = W^{[n+1]} \left( f_n \circ f_{n-1} \circ \cdots \circ f_1 (W^{[0]}x + b^{[0]}) \right) + b^{[n+1]},
\end{equation}
where $\theta$ denote the sets of parameters, {\it i.e.}, 
\begin{equation*}
	\theta = \left\{W^{[0]}, b^{[0]}, W^{[1_1]}, b^{[1_1]}, W^{[1_2]}, b^{[1_2]} , \cdots
	, W^{[n_1]}, b^{[n_1]}, W^{[n_2]}, b^{[n_2]} ,  W^{[n+1]},b^{[n+1]} \right\}.
\end{equation*}
In Figure \ref{fig:resnet}, we plot the architecture of  the ResNet with 4 block. 
 
\subsection{Deep unfitted Nitsche formulation}
The main advantage of the unfitted Nitsche's method lies in the ability to use meshes independent of the location of the interface. This is possible by employing two different ansatz functions on the interface elements (elements cut through by the interface): one is for the interior domain, and the other one is for the exterior domain. Those two different ansatz functions are discontinuous across the interface $\Gamma$ and patched together by Nitsche's method. In the traditional unfitted Nitsche's methods, the ansatz functions are piecewise polynomials. Following this line, we adopt two different deep neural network functions as the two ansatz functions to minimize the unfitted Nitsche energy \eqref{equ:energy}.

   Let  $u_{\theta_i}$ be the ansatz function in $\Omega_i$ ($i=1,2$) and denote $u_{\Theta}(x)=(u_{\theta_1}(x), u_{\theta_2}(x))$.
 where $\Theta=(\theta_1, \theta_2)$.  Before we proceed, we should make sure $u_{\Theta}$ satisfies the Dirichlet boundary condition 
  \ref{equ:bnd}.  For such purpose, we introduce an additional boundary penalty term as 
  \begin{equation}
  	L_b(u_{\Theta}) = \int_{ \partial  \Omega} |u_{\Theta}-g|^2ds. 
  \end{equation}
  Ideally, we expect $L_b$ close to zero.  We define the extended unfitted Nitsche's functional    $L$ as 
  \begin{equation}
  	L(u_{\Theta}) = L_n(u_{\Theta}) + \gamma_b L_b(u_{\Theta}),
  \end{equation}
  where $\gamma_b$ is the boundary penalty parameter. 
  
   Then, our deep unfitted Nitsche method is equivalent the following optimization problem 
 \begin{equation} \label{equ:optimize}
 	\min_{\Theta} L(u_{\Theta}).
 \end{equation}

 \begin{remark}
 	 In the literature, there are several alternative methods to impose the Dirichlet boundary conditions by building the Dirichlet boundary condition into the loss function \cite{LiMi2021} and training another deep neural network \cite{ShYa2021}.  For the sake of simplicity, we choose to impose the Dirichlet boundary condition for deep neural network functions using the penalty method as in \cite{EYu2018}. 
 	 Interesting readers are referred to \cite{ChDW2020} for a comparison study of different boundary conditions handling methods.   
 \end{remark}

To solve the optimization problem using stochastic gradient descent type algorithms ({\it e.g.}, SGD\cite{GoBC2016} or  ADAM\cite{KiBa2015adam}), we approximate the integrals by using the Monte-Carlo method, where the number of integral point is independent of the dimensional of the underlying domain.  Suppose $\{x_k^{i}\}_{k=1}^{N_i}$ are the uniformly sampled point in the domain $\Omega_i$ for $i=1,2$.  Similarly, let $\{x_k^{f}\}_{k=1}^{N_f}$ and $\{x_k^{b}\}_{k=1}^{N_b}$ be the randomly sampled point on the interface $\Gamma$ and the domain boundary $\partial \Omega$, respectively.   The loss function $\hat{L}$ is defined as 
\begin{equation}\label{equ:loss}
	\begin{split}
		\hat{L}(u_{\Theta}) = & \frac{|\Omega_1|}{N_1}\sum_{k=1}^{N_1}
		\left(\frac{\beta_1}{2}\nabla u_{\theta_1}(x_k^1) \cdot 
		\nabla u_{\theta_1}(x_k^1) - f(x_k^1) u_{\theta_1}(x_k^1)\right) + \\
		&\frac{|\Omega_2|}{N_2}\sum_{k=1}^{N_2}
		\left(\frac{\beta_2}{2}\nabla u_{\theta_2}(x_k^2) \cdot 
		\nabla u_{\theta_2}(x_k^2) - f(x_k^2) u_{\theta_2} (x_k^2)\right) + \\
		& \frac{|\Gamma|}{N_f}\sum_{k=1}^{N_f}\left( \frac{\gamma_f}{2}(\llbracket u_{\Theta}(x_k^f)\rrbracket - p(x_k^f))^2
		 - \frac{\gamma_f}{2} p(x_k^f)^2 \right) + \\
		 & \frac{|\Gamma|}{N_f}\sum_{k=1}^{N_f}\left( (p(x_k^f) - \llbracket u_{\Theta}(x_k^f)\rrbracket)
 	 \dgal{\beta \partial_n u_{\Theta}(x_k^f)}\right) -\\
 	 & \frac{|\Gamma|}{N_f}\sum_{k=1}^{N_f}\left( q(x_k^f)\dgal{\beta \partial_n u_{\Theta}(x_k^f)}^{\ast}  \right) + \\
 	 & \frac{|\partial\Omega|}{N_b}\sum_{k=1}^{N_b} \gamma_b\left( u_{\theta_2}(x_k^b) - g(x_k^b)\right)^2,
	\end{split}
\end{equation}
 where $|\Omega_i|$ is the measure of $\Omega_i$ in $\mathbb{R}^d$ ($i=1,2$)  and $|\Gamma|$ (or $|\partial\Omega|$)  is the measure of $\Gamma$ (or $\partial \Omega$) in $\mathbb{R}^{d-1}$.  Then, the discrete counterpart  of the optimization problem  \eqref{equ:optimize} reads as 
\begin{equation}\label{equ:disopt}
	\min_{\Theta} \hat{L}(u_{\Theta}).
\end{equation}

The discrete optimization problem \eqref{equ:disopt} actually involves the training of two deep neural network functions $u_{\theta_1}$ and $u_{\theta_2}$. Those two neural networks can be trained independently using the same loss function $\hat{L}$.  The gradient of deep neural network function can be efficiently calculated using automatic differentiation functionality  \cite{pytorchad2019} in the modern machine learning platform.
 
 In the loss function\eqref{equ:loss},  we need to compute the measure of each $\Omega_i$. For problems with simple geometric shapes, we can compute them analytically.  In general, we can Monte-Carlo simulation like the hit-or-miss method to estimate the measure. Similarly, we can estimate the measure of $\Gamma$ and $\partial\Omega$.

 \section{Numerical Experiments}\label{sec:example}
 In this section, we present several numerical examples to illustrate the performance of the proposed deep unfitted Nitsche method. 
The proposed algorithm is implemented using the machine learning platform Pytorch \cite{NEURIPS2019_9015}. In the following numerical experiments, we choose $u_{\theta_1}$ and $u_{\theta_2}$ have the same neural network architectures and select the activation function $\sigma=\tanh$. Both deep neural networks are initialized by Xavier initialization to prevent from exploding or vanishing and are trained independently using ADAM\cite{KiBa2015adam}.  In all the following numerical experiments, we choose the learning rate $lr=0.001$ and generate new mini-batches every 10 epochs.  To provide a qualitative description of training results, we  use  the following relative $L^2$-error
\begin{equation}
	\text{Error}  = \frac{\|u-u_{\Theta}\|_{0, \Omega_1\cup \Omega_2}}{\|u\|_{0, \Omega_1\cup \Omega_2}}.
\end{equation}
 The relative $L^2$-error is computed using Monte-Carlo methods with $10,000$ uniform sampled points.  The relative $L^2$ error are computed and recorded in every 100 epochs. Similarly, we record the loss in every 100 epochs. 
 
 \subsection{Flower shape interface problem in 2D} 
 In our first example, we consider the flower interface problem in the domain $\Omega=(-1,1)\times(-1,1)$ as in \cite{GuYa2018, guo2017gradient}.  The flower shape interface $\Gamma$ is given by  the following polar coordinate 
 \begin{equation}\label{equ:flower}
 	r = \frac{1}{2} + \frac{\sin(5\theta)}{7}. 
 \end{equation}
The piecewise diffusion coefficient are  $\beta_1=1$ and $\beta_2=10$. We choose the right hand side function to fit the exact solution 
\begin{equation*}
u(x) =
\left\{
\begin{array}{ll}
    e^{(x_1^2+x_2^2)}, & \text{if } x=(x_1, x_2)\in \Omega_1,\\
    0.1(x_1^2+x_2^2)^2-0.01\ln(2\sqrt{x_1^2+x_2^2}),& 
     \text{if } x=(x_1, x_2)\in \Omega_2.\\
   \end{array}
\right.
\end{equation*}
 The nonhomogeneous jump conditions \eqref{equ:fluxjump} and \eqref{equ:valuejump} can be calculated from the above exact solution.

 \begin{figure}[!h]
   \centering
  \subcaptionbox{\label{fig:flower_sol_nn}}
   {\includegraphics[width=0.47\textwidth]{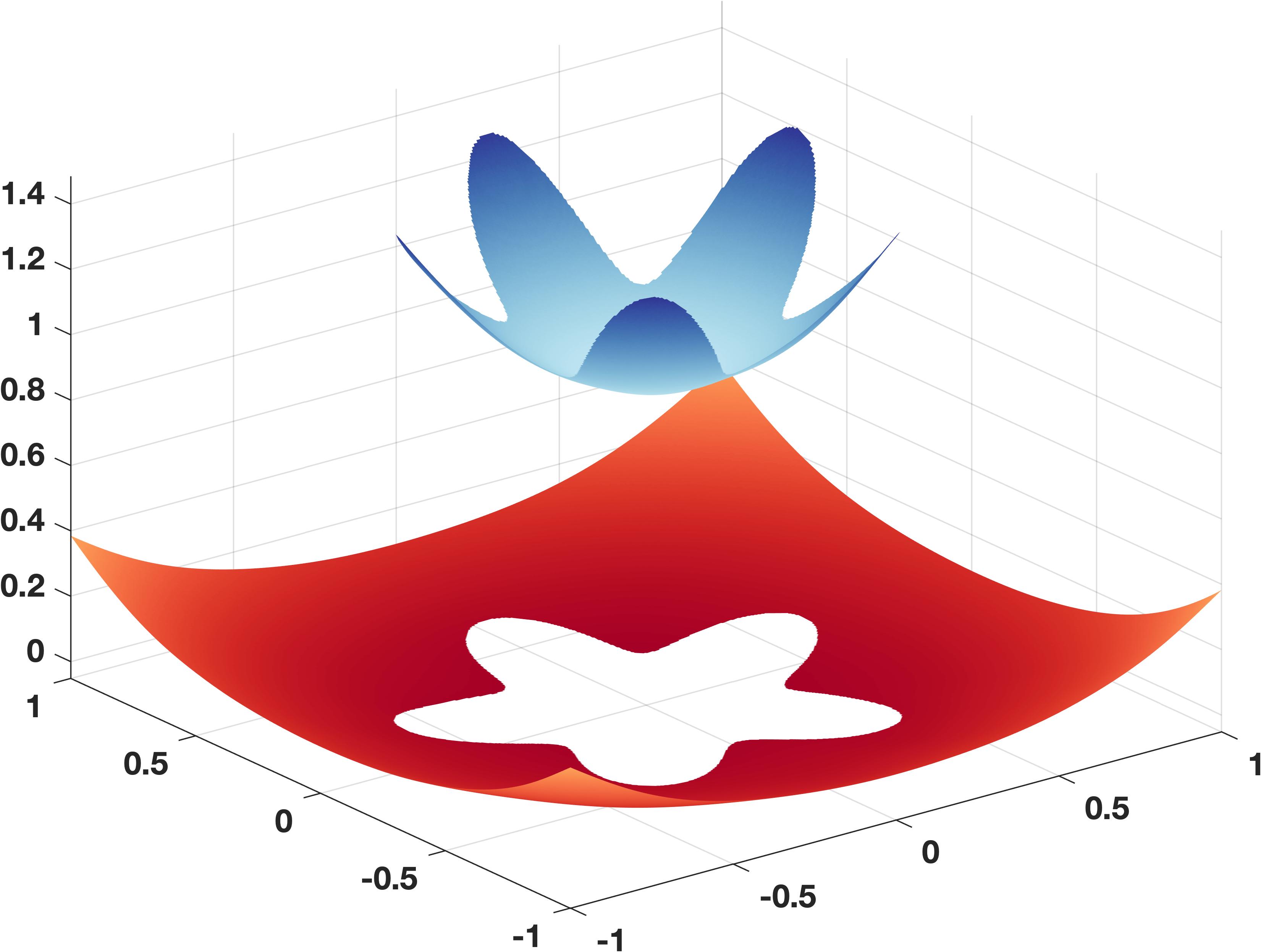}}
   \subcaptionbox{\label{fig:fower_sol_nn}}
   {\includegraphics[width=0.47\textwidth]{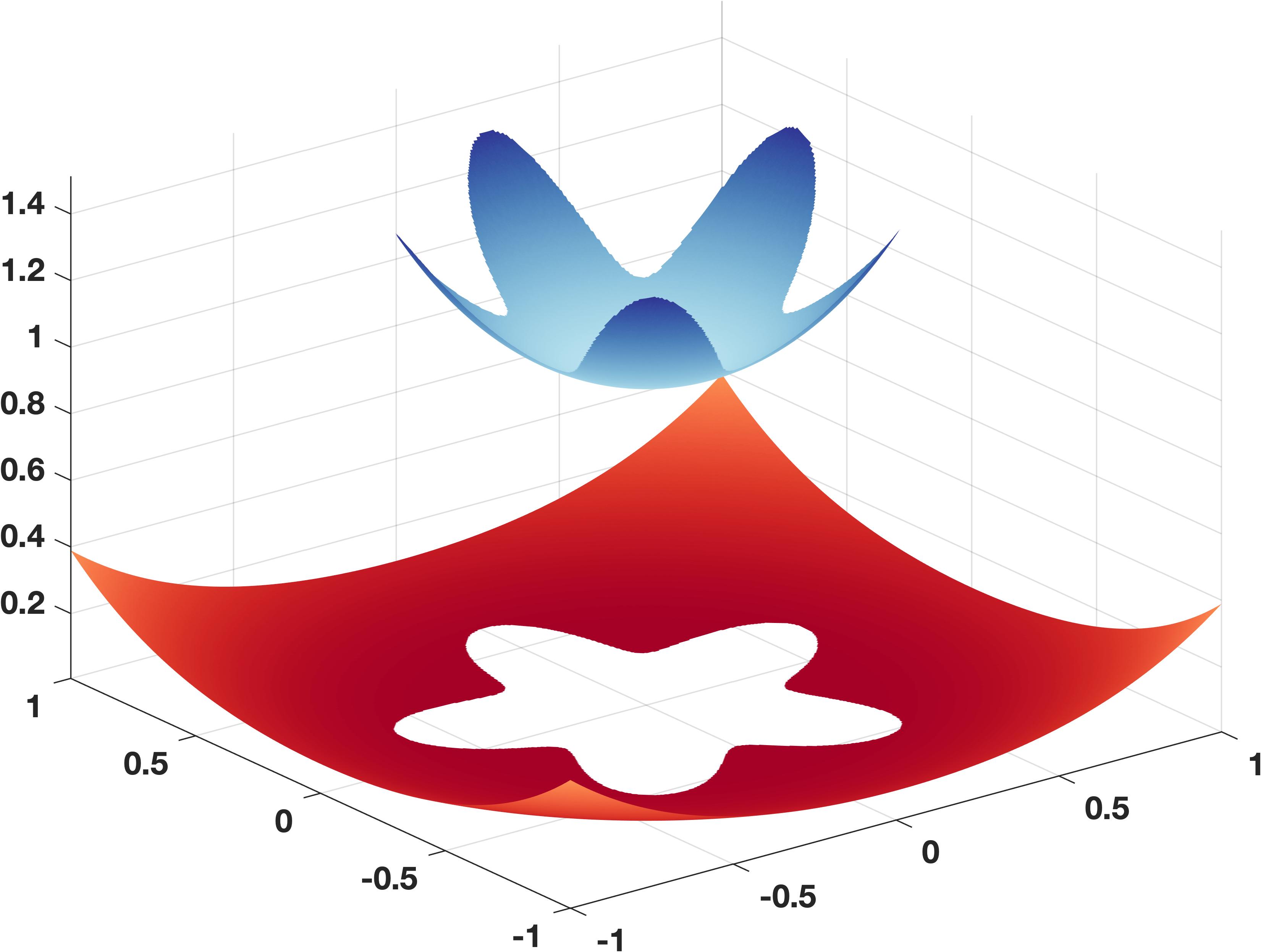}}
   \caption{Comparison of solutions for flower interface problem: (a) Deep unfitted Nitsche's solution; (b) Exact solution. Although there is an inhomogeneous jump of funciton values at the interface, both solutions match well.}\label{fig:flower_sol}
\end{figure}

  \begin{figure}[!h]
   \centering
  \subcaptionbox{\label{fig:flower_error}}
   {\includegraphics[width=0.47\textwidth]{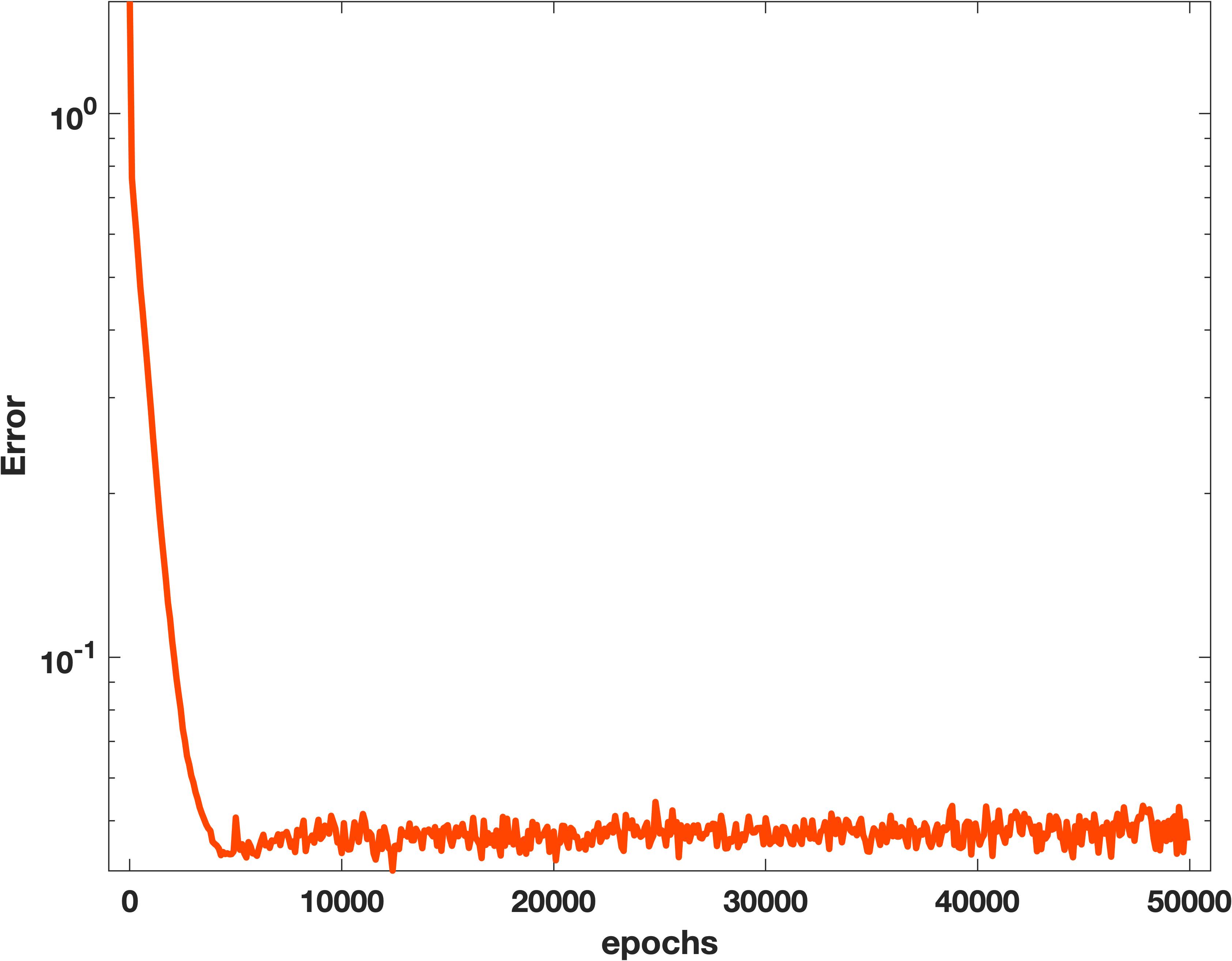}}
   \subcaptionbox{\label{fig:fower_loss}}
   {\includegraphics[width=0.47\textwidth]{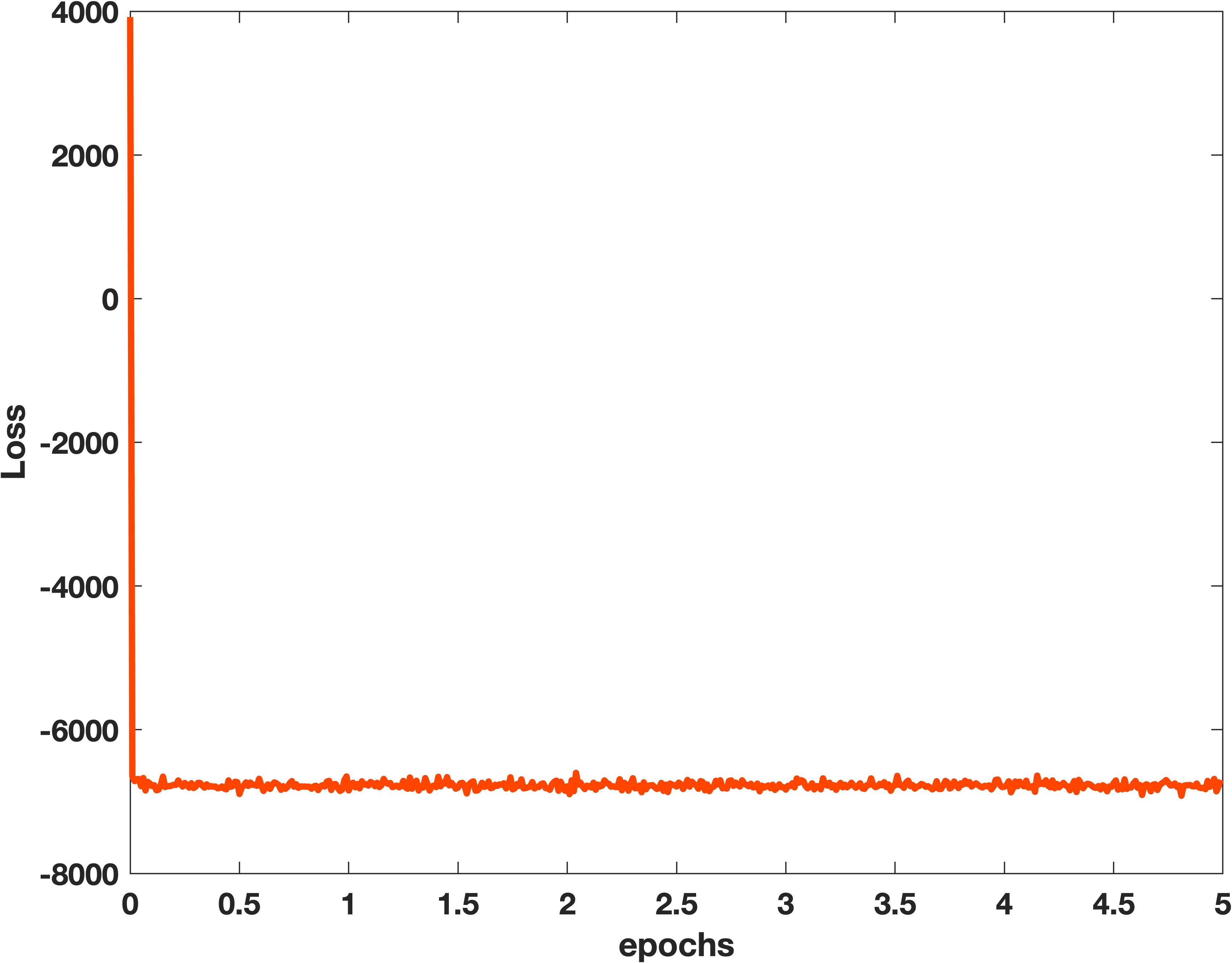}}
   \caption{Training process of flower interface problem: (a) Decay of the relative $L^2$-error; (b) Decay of the loss. }\label{fig:flower_sum}
\end{figure}

 In this case, one can compute $|\Omega_1|= \frac{51}{192}\pi$ and $\Omega_2 = 4 - \frac{51}{192}\pi$. To general uniformly distributed random points in $\Omega_1$ and $\Omega_2$, we firstly generate uniformly distributed random points in the whole domain $\Omega$. Then, we count the random points inside the interface $\Gamma$ as the randomly sampled points in $\Omega_1$ and the rest random points are in $\Omega_2$. The randomly sampled points $\Gamma$ is produced by generating the uniformed sampled points on the interval $(0, 2\pi)$ and then are mapped onto the interface $\Gamma$ using \eqref{equ:flower}.

In this numerical test,  we choose the  ResNet with 3 blocks with $m=10$ for each ansatz function, as illustrated in the diagram \ref{fig:resnet}.  Each ansatz function has 701 parameters. We uniformly sample 1024 points in the domain $\Omega$. It turns out that there are $219$ points inside the interface $\Gamma$. In other words, $N_1=219$ and $N_2=805$. In addition, we just choose $N_f = 256$ and $N_b = 128$. We choose $\gamma_b = 5000$ and $\gamma_f=1000$. The corresponding decay of errors and loss functions during the training process are plotted in Figure \ref{fig:flower_sum}. We can see the error decays very fast at the initial several thousand epochs and then fluctuate. After 50000 epochs, the relative $L^2$-error is reduced to $4.6\%$. In Figure \ref{fig:flower_sol}, we present a comparison between the deep unfitted Nitsche method(DUNM) solution and the exact solution. We can see the DUNM solution match well with the exact solution even though there is an inhomogeneous jump in the function values.

  \begin{figure}[!h]
   \centering
  \subcaptionbox{\label{fig:cp1000m1_error}}
   {\includegraphics[width=0.47\textwidth]{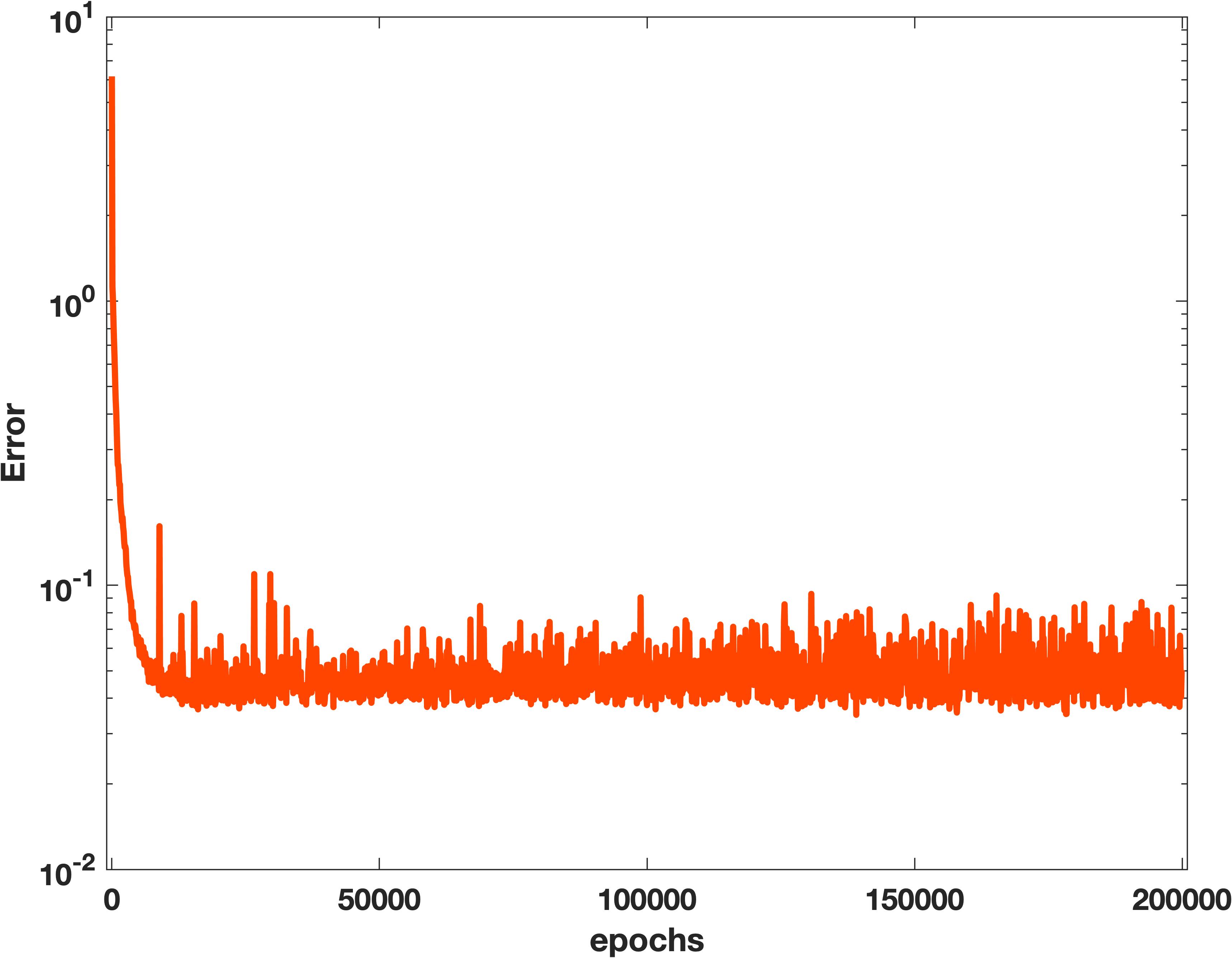}}
   \subcaptionbox{\label{fig:cp1000m1_loss}}
   {\includegraphics[width=0.47\textwidth]{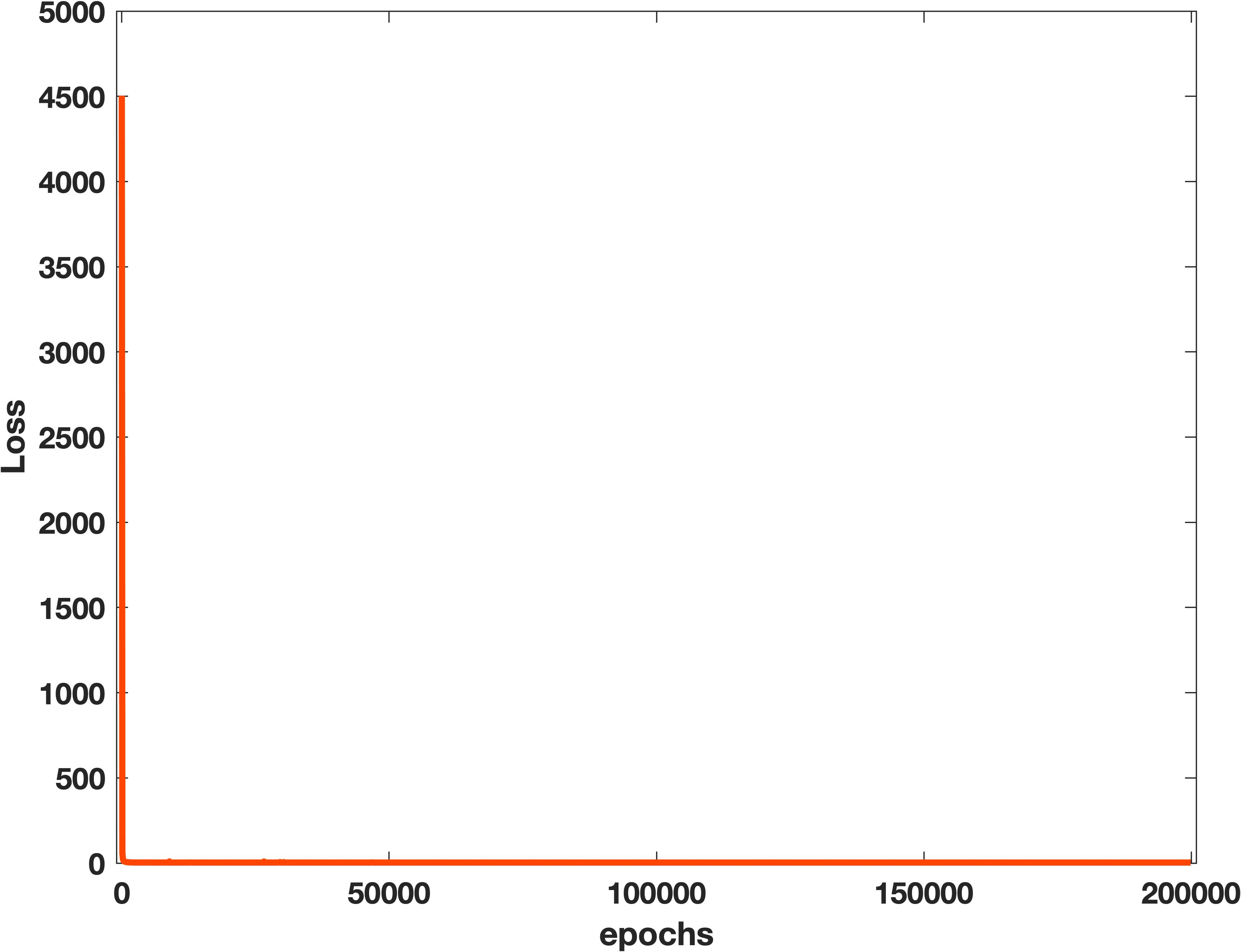}}
   \caption{Training process of circle interface problem with $\beta_1 = 1$ and $\beta_2 = 1000$: (a) Decay of the relative $L^2$-error; (b) Decay of the loss. }\label{fig:cp1000m1}
\end{figure}

  \begin{figure}[!h]
   \centering
  \subcaptionbox{\label{fig:cp100000m1_error}}
   {\includegraphics[width=0.47\textwidth]{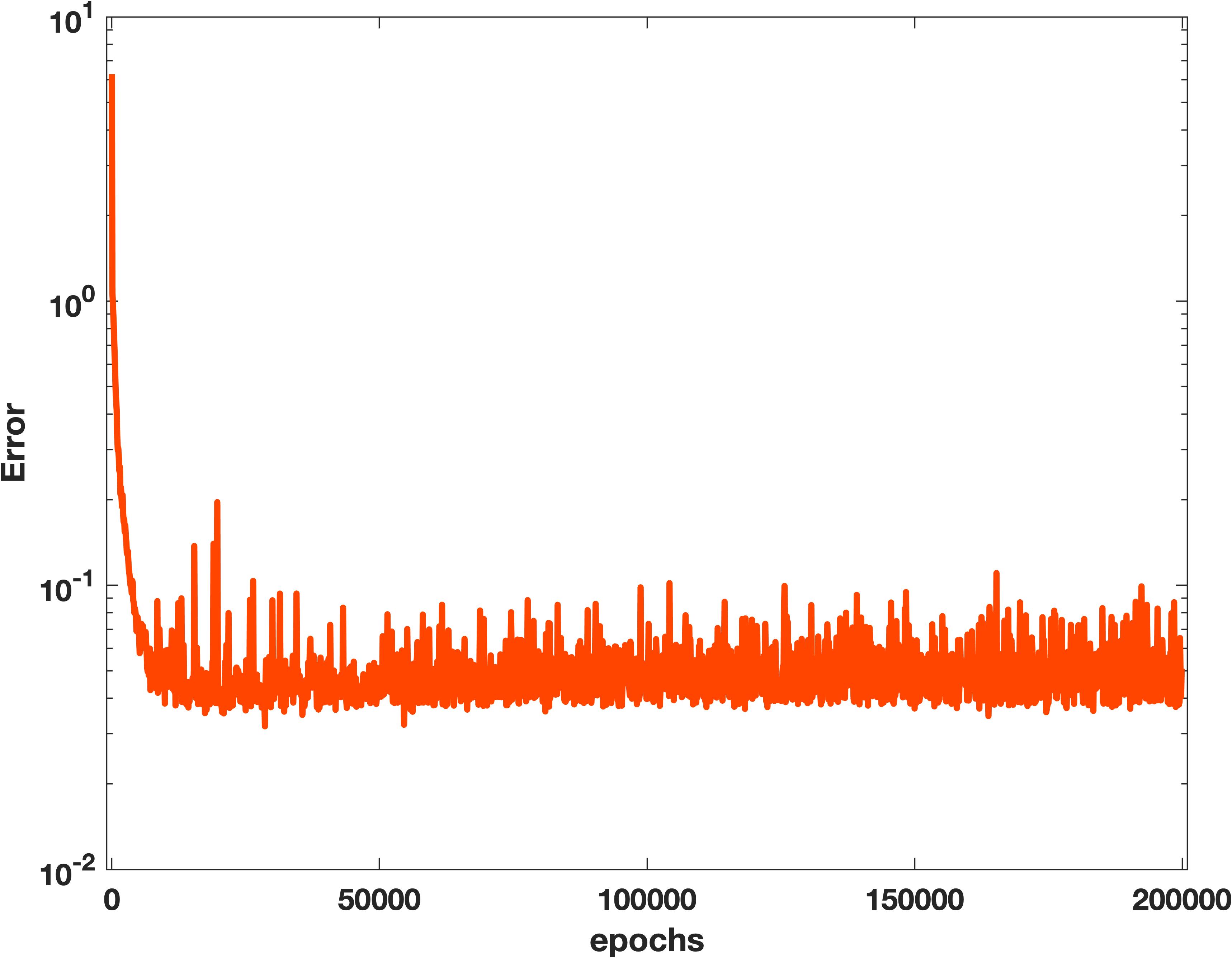}}
   \subcaptionbox{\label{fig:cp100000m1_loss}}
   {\includegraphics[width=0.47\textwidth]{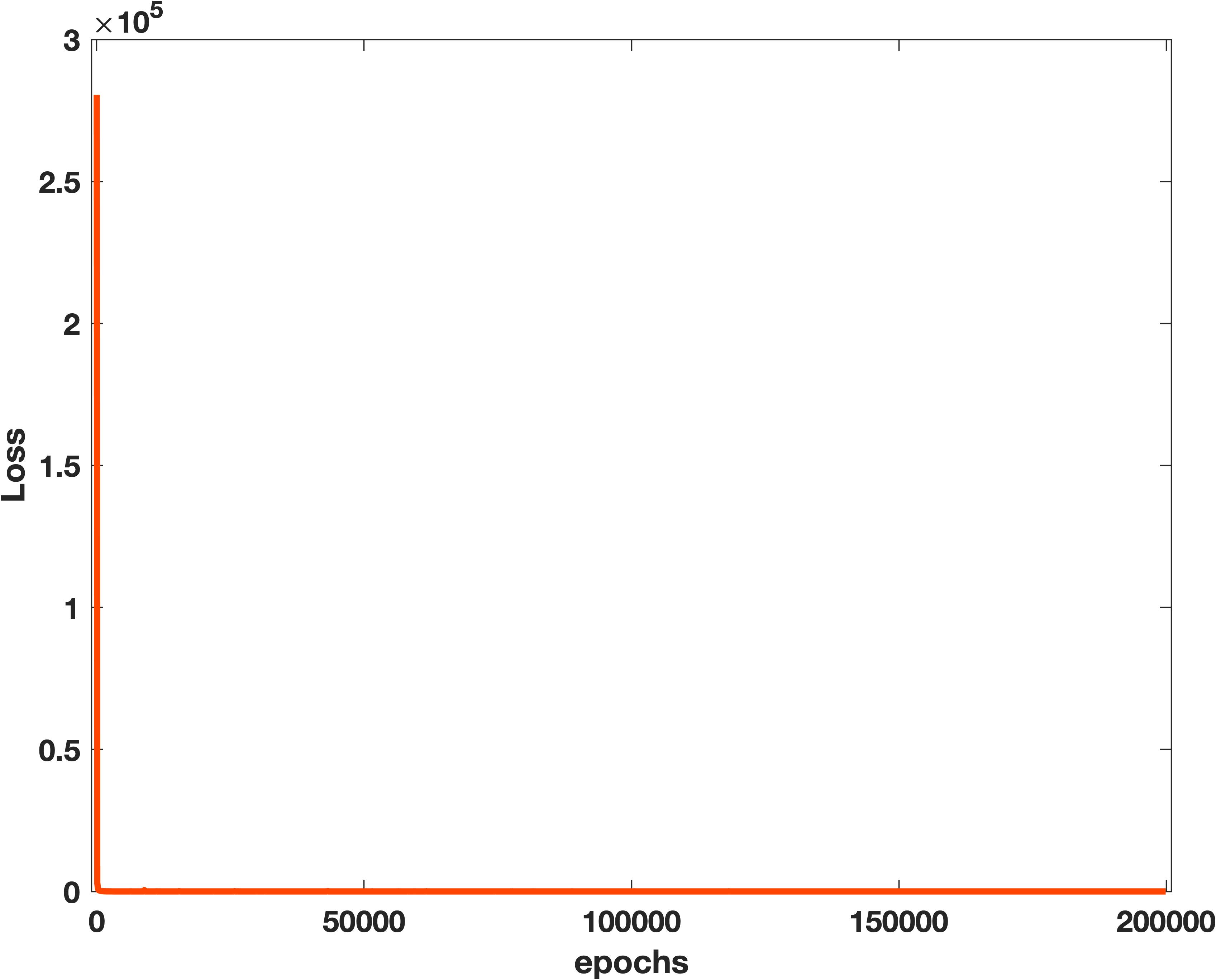}}
   \caption{Training process of Training process of circle interface problem with $\beta_1 = 1$ and $\beta_2 = 100000$: (a) Decay of the relative $L^2$-error; (b) Decay of the loss. }\label{fig:cp100000m1}
\end{figure}

  \begin{figure}[!h]
   \centering
  \subcaptionbox{\label{fig:cp1m1000_error}}
   {\includegraphics[width=0.47\textwidth]{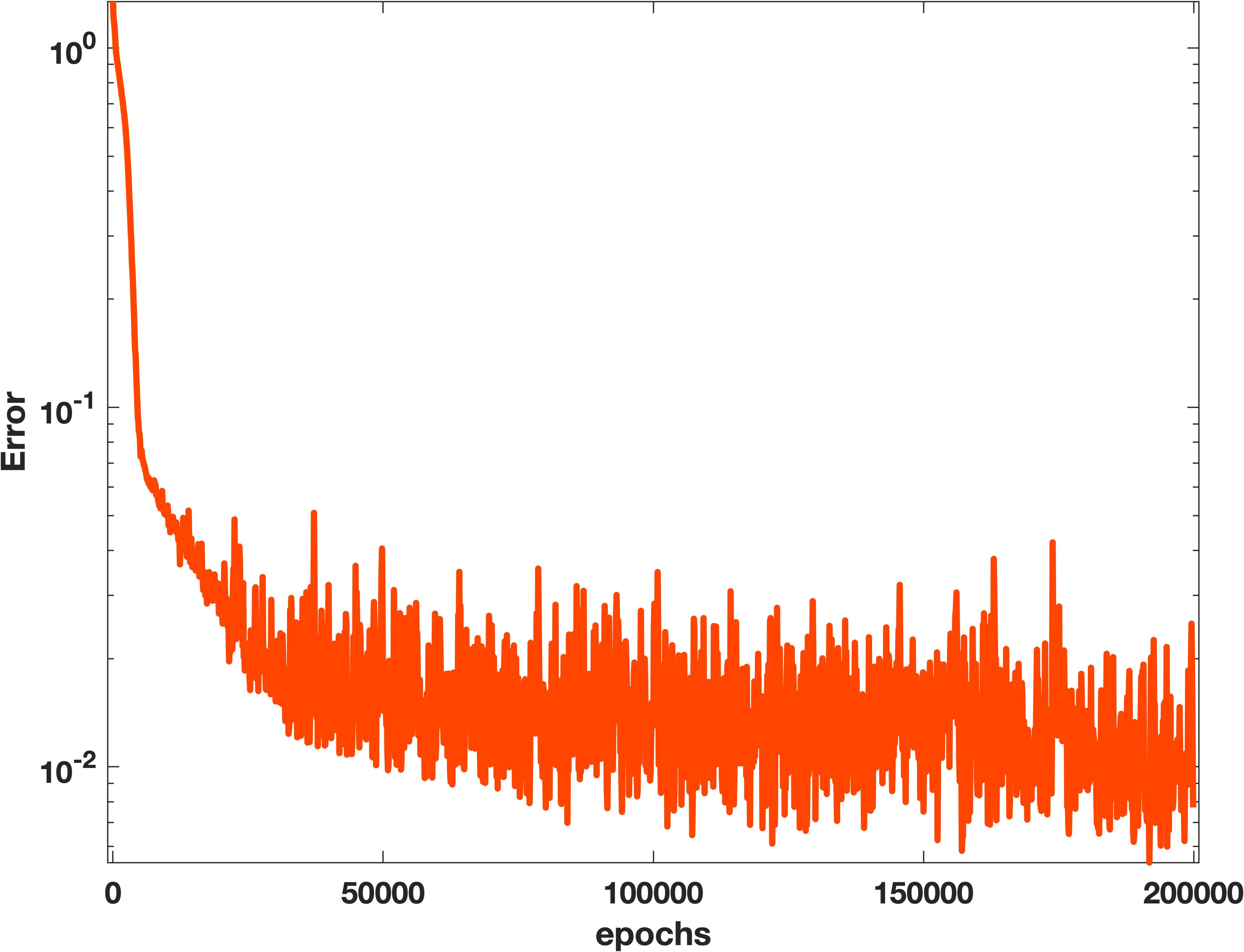}}
   \subcaptionbox{\label{fig:cp1m1000_loss}}
   {\includegraphics[width=0.47\textwidth]{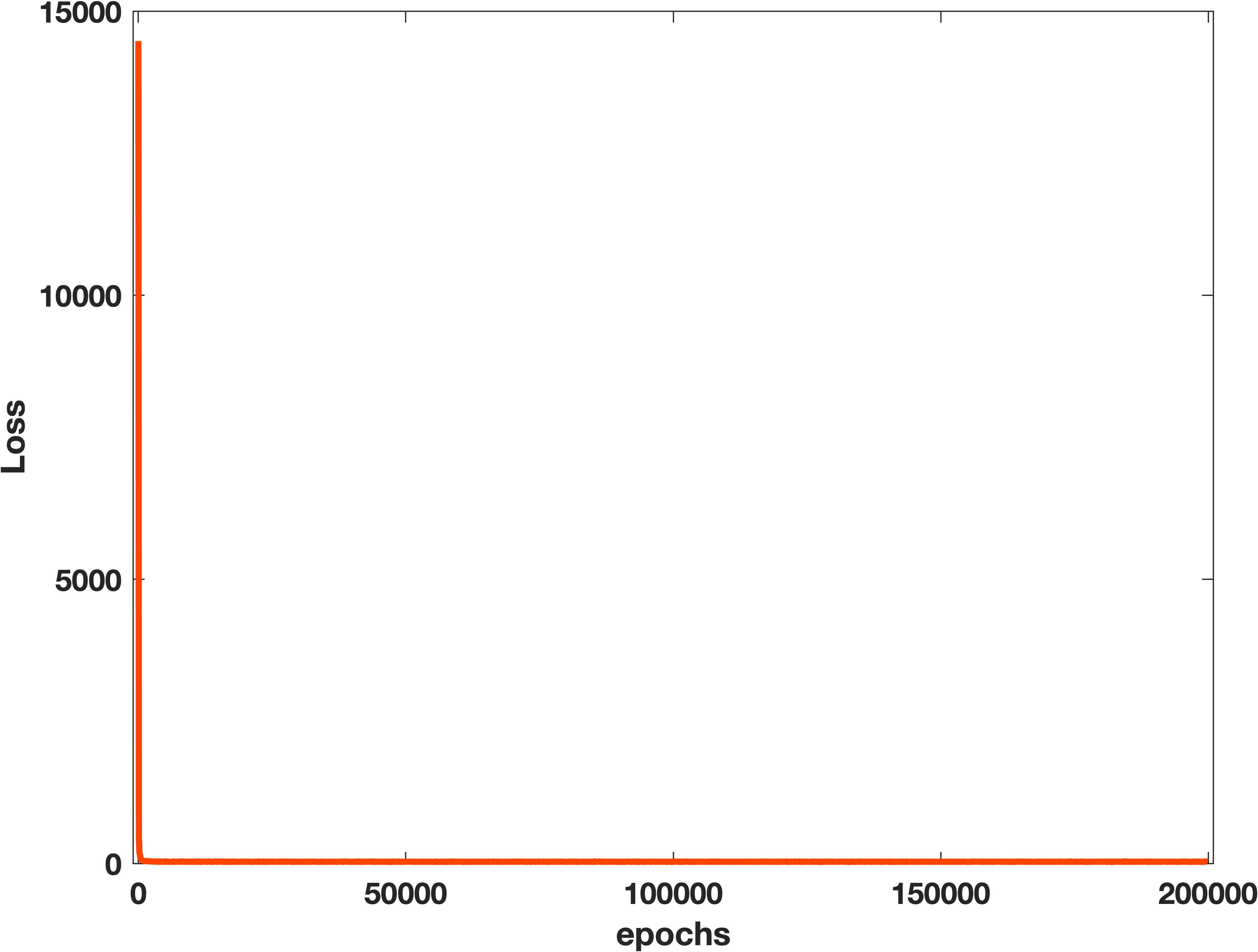}}
   \caption{Training process of circle interface problem with $\beta_1 = 1000$ and $\beta_2 = 1$: (a) Decay of the relative $L^2$-error; (b) Decay of the loss. }\label{fig:cp1m1000}
\end{figure}

  \begin{figure}[!h]
   \centering
  \subcaptionbox{\label{fig:cp1m100000_error}}
   {\includegraphics[width=0.47\textwidth]{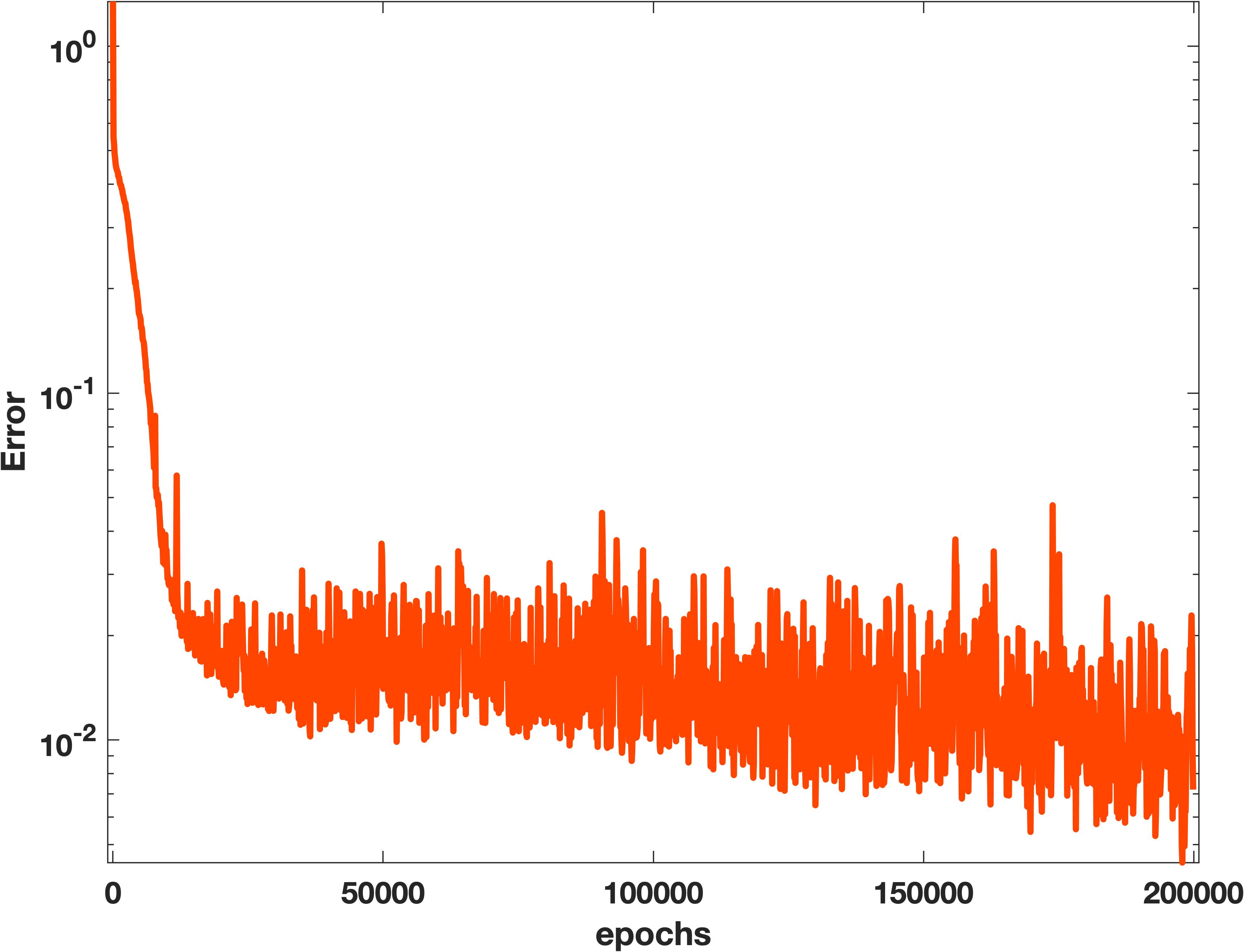}}
   \subcaptionbox{\label{fig:cp1m100000_loss}}
   {\includegraphics[width=0.47\textwidth]{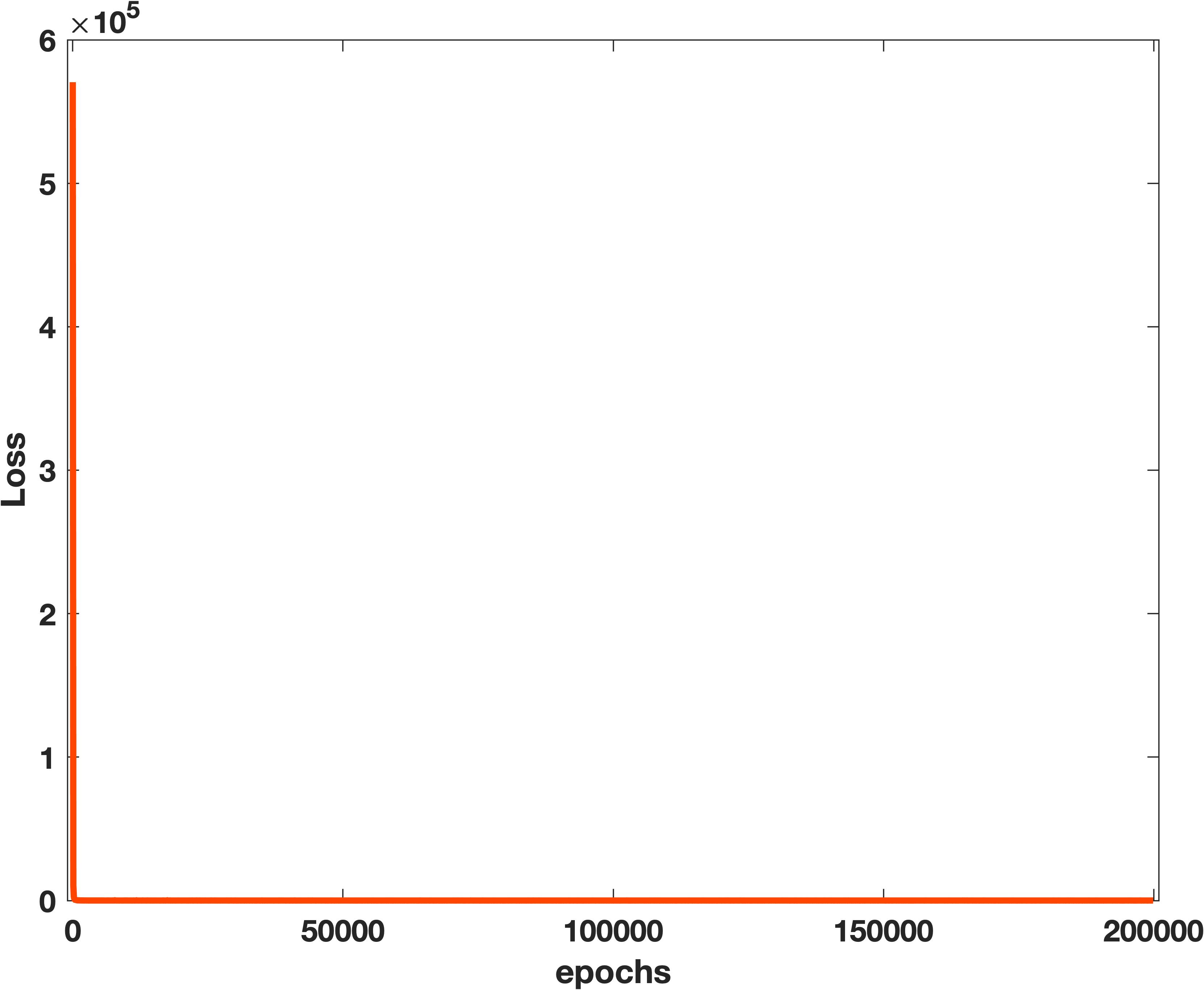}}
   \caption{Training process of circle interface problem with $\beta_1 = 100000$ and $\beta_2 = 1$: (a) Decay of the relative $L^2$-error; (b) Decay of the loss. }\label{fig:cp1m100000}
\end{figure}

 \begin{figure}[!h]
   \centering
  \subcaptionbox{\label{fig:cp1000m1_sol_nn}}
   {\includegraphics[width=0.47\textwidth]{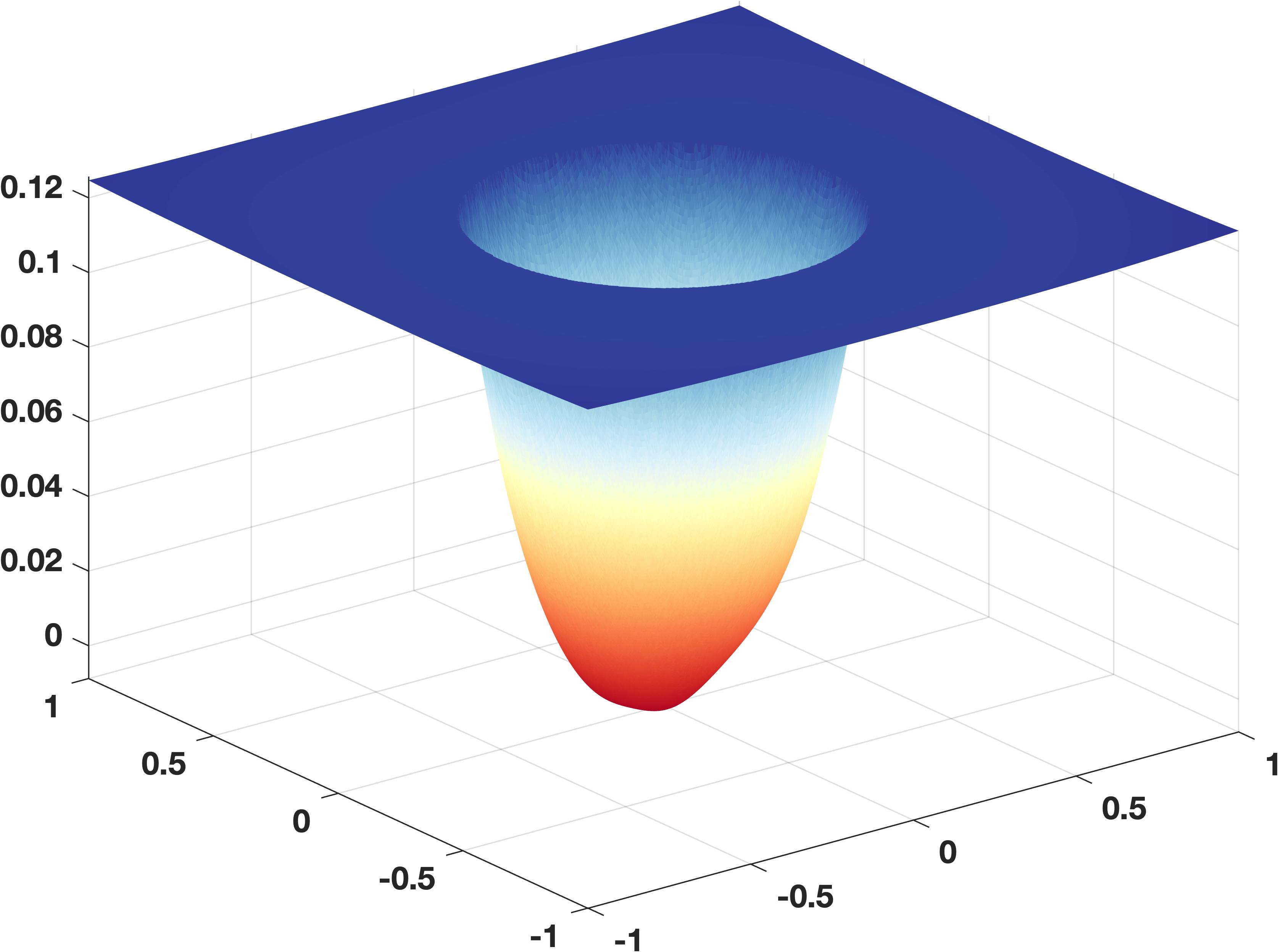}}
   \subcaptionbox{\label{fig:cp1000m1_sol_exact}}
   {\includegraphics[width=0.47\textwidth]{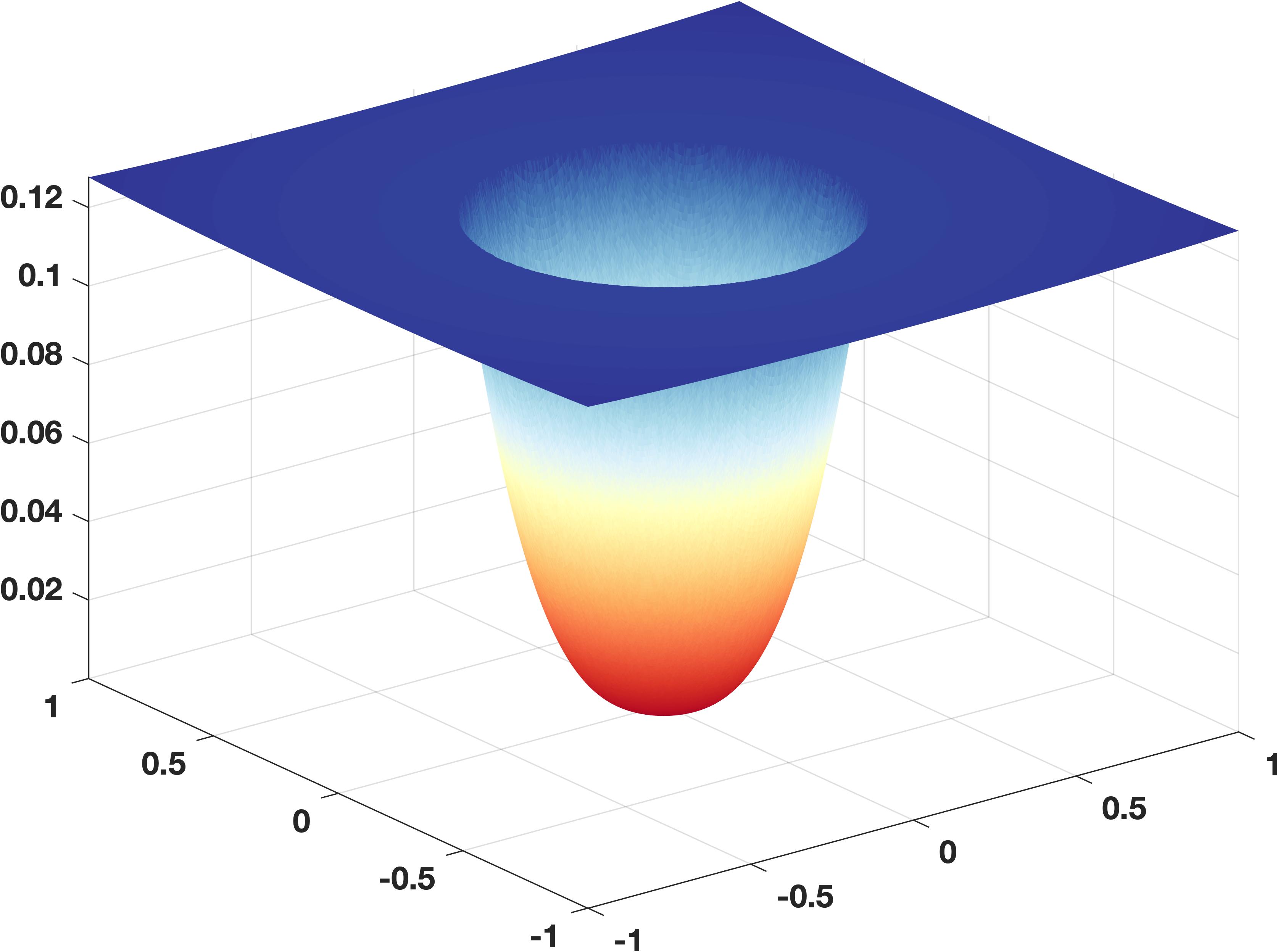}}
   \caption{Comparison of solutions for circular interface problem with $\beta_1=1$ and $\beta_2=1000$: (a) Deep unfitted Nitsche's solution; (b) Exact solution. The solutions match well for this high-contrast interface example. }\label{fig:circlep1000_m1_sol}
\end{figure}

 \begin{figure}[!h]
   \centering
  \subcaptionbox{\label{fig:cp1m1000_sol_nn}}
   {\includegraphics[width=0.47\textwidth]{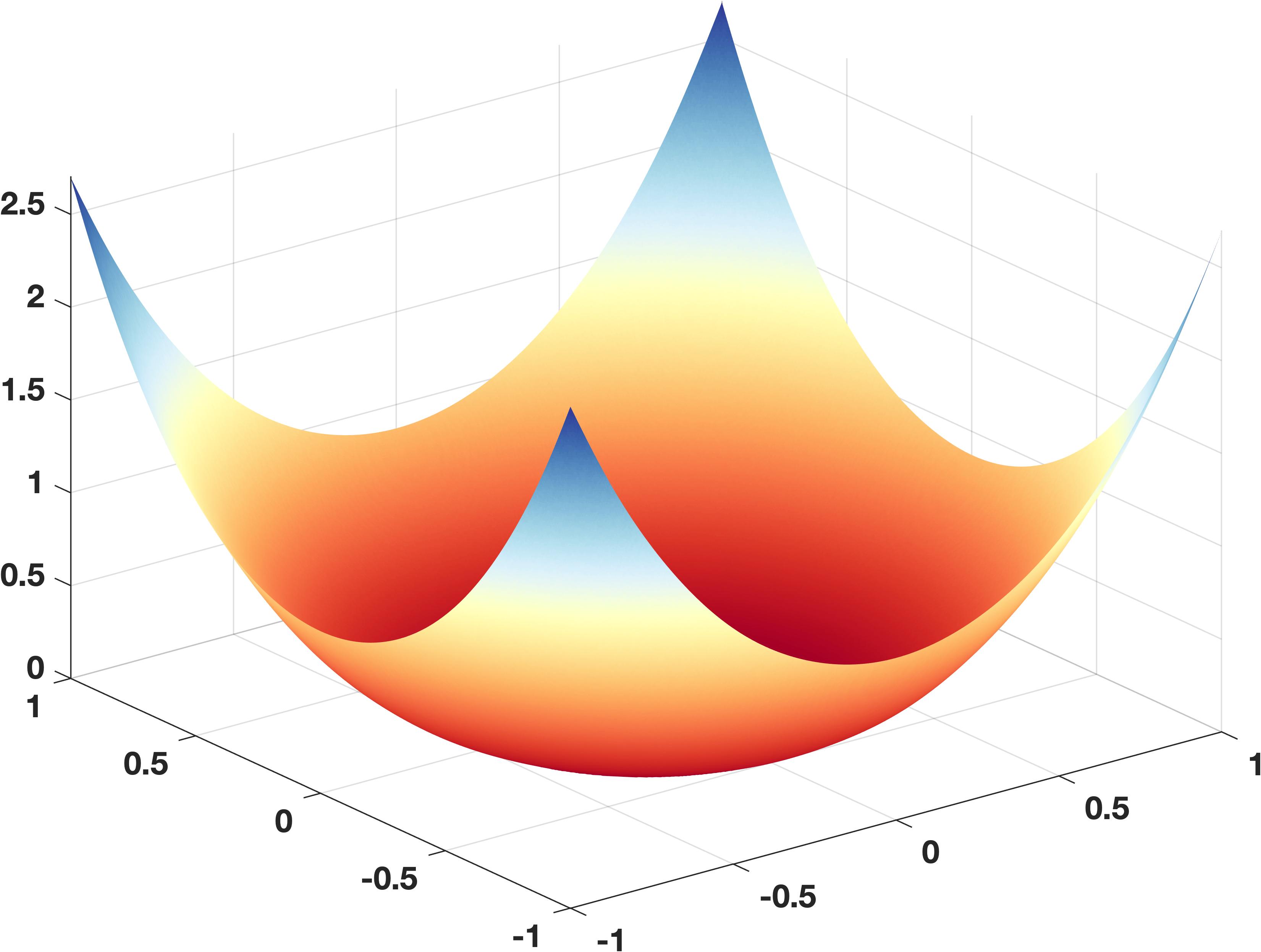}}
   \subcaptionbox{\label{fig:cp1m1000_sol_exact}}
   {\includegraphics[width=0.47\textwidth]{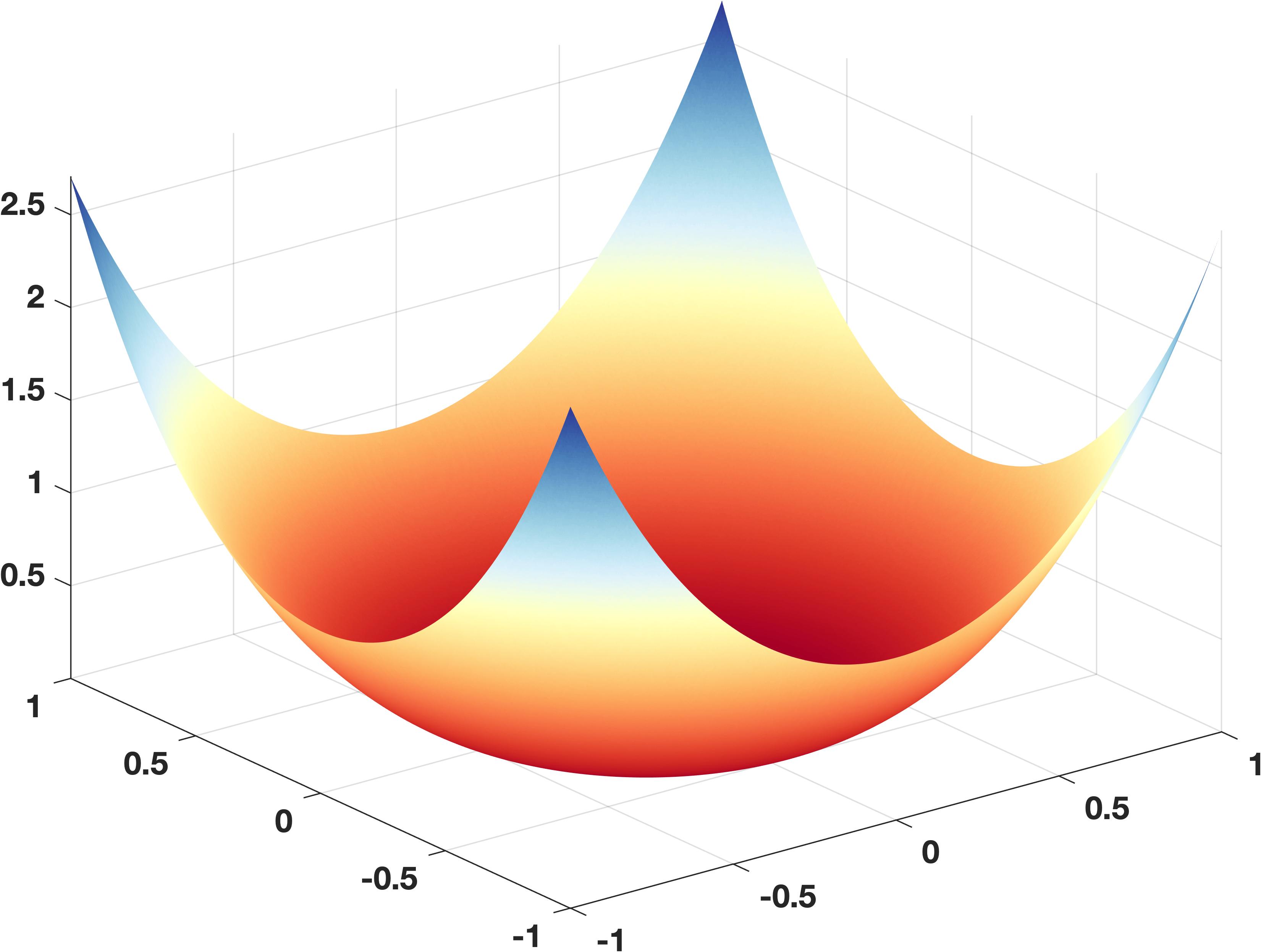}}
   \caption{Comparison of solutions for circular interface problem with $\beta_1=10000$ and $\beta_2=1$: (a) Deep unfitted Nitsche's solution; (b) Exact solution. }\label{fig:circlep1_m1000_sol}
\end{figure}

 \subsection{High-contrast interface problem in 2D} 
 In this example, we consider the high contrast interface problem with homogeneous jump conditions in the domain $\Omega=(-1,1)\times(-1,1)$ as in  \cite{GuYa2018, guo2018gradient, guo2017gradient}. The interface $\Gamma$ is the circle of radius $r_0$ centered at the original.  The exact solution is 
 \begin{equation*}
u(x_1,x_2) =
\left\{
\begin{array}{ll}
    \frac{r^3}{\beta_1}   &  \text{if }   x\in \Omega_1, \\
      \frac{r^3}{\beta_2} + \left( \frac{1}{\beta_1}-\frac{1}{\beta_2} \right)r_0^3&  \text{if }  x\in \Omega_2,\\
   \end{array}
\right.
\end{equation*}
where $r = \sqrt{x_1^2 + x_2^2}$.

In this example, we consider the  circle with radius $r_0=0.5$. We also use the ResNet with 3 blocks with $10$ output in each block to represent each function. 
The uniformly random sampled points are generated by the same method as in previous example. Again, we generate $1024$ random sampled points in $\Omega$. In that case, $N_1=212$ points are inside $\Gamma$ and $N_2=812$ points are outside $\Gamma$. Similarly,  $N_f = 256$ and $N_b = 128$. The penalty parameters $\gamma_b = 5000$ and $\gamma_f=1000$ are the exactly the same as previous case.

\begin{table}[h!] 
\centering
 \caption{Relative $L^2$ errors for high contrast interface problems}\label{tab:highcontrast}
 \begin{tabular}{|c | c | c | c | c|} 
 \hline
 $\beta_1/\beta_2$ & $1000/1$ & $1/1000$ & $100000/1$ & $1/100000$ \\  \hline
Error($\%$) & 2.29 & 0.62 & 2.36 & 0.63 \\
 \hline
 \end{tabular}
\end{table}

In the following numerical tests, we focus on the training  of high contrast interface problem by considering the following  four typical different
jump ratios:  $\beta_1/\beta_2 = 1000/1$ (large jump), $\beta_1/\beta_2 = 1/1000$ (large jump), $\beta_1/\beta_2 = 1/100000$ (huge jump),
and $\beta_1/\beta_2 = 100000$ (huge jump). The training processes are described in Figures \ref{fig:cp1000m1}-\ref{fig:cp1m1000} for each cases, respectively. In these figures, we notice that the errors quickly decay to some specific errors and then oscillate around them.  The relative $L^2$ errors after 200000 epochs are summarized in the Table~\ref{tab:highcontrast}, where one can observe that the errors are not influenced by the jump ratios.

In the Figures~\ref{fig:circlep1000_m1_sol} and  \ref{fig:circlep1_m1000_sol}, we plot the DUNM solution and the exact solution for the case $\beta_1/\beta_2=1000/1$ and $\beta_1/\beta_2=1/1000$, respectively.  It is clear to see that the DUNE solutions match well with the corresponding exact solutions.  We have also compared two other cases and observed the same phenomena.

 \subsection{High-dimensional interface problems}
In this example, we consider the $d$-dimensional interface problem in the unit cube  $\Omega = [-0.5, 0.5]^d$.  The unit cube is divided  into two parts $\Omega_1$ and $\Omega_2$ by the d-dimensional sphere with radius $r_0 = 0.4$. The diffusion coefficients are $\beta_1=1$ and $\beta_2=10$.  The exact solution  is 
 \begin{equation*}
u(x) =
\left\{
\begin{array}{ll}
    r^3   &  \text{if }   x\in \Omega_1, \\
      \frac{r^3}{10} + 0.0576&  \text{if }  x\in \Omega_2,\\
   \end{array}
\right.
\end{equation*}
where $r = |x| = \sqrt{\sum_{j=1}^{d}x_j^2}$.  Therefore, we have homogeneous jump conditions. The right hand side function and boundary condition can be obtained from $u$.

  \begin{figure}[!h]
   \centering
  \subcaptionbox{\label{fig:sphere3d_error}}
   {\includegraphics[width=0.47\textwidth]{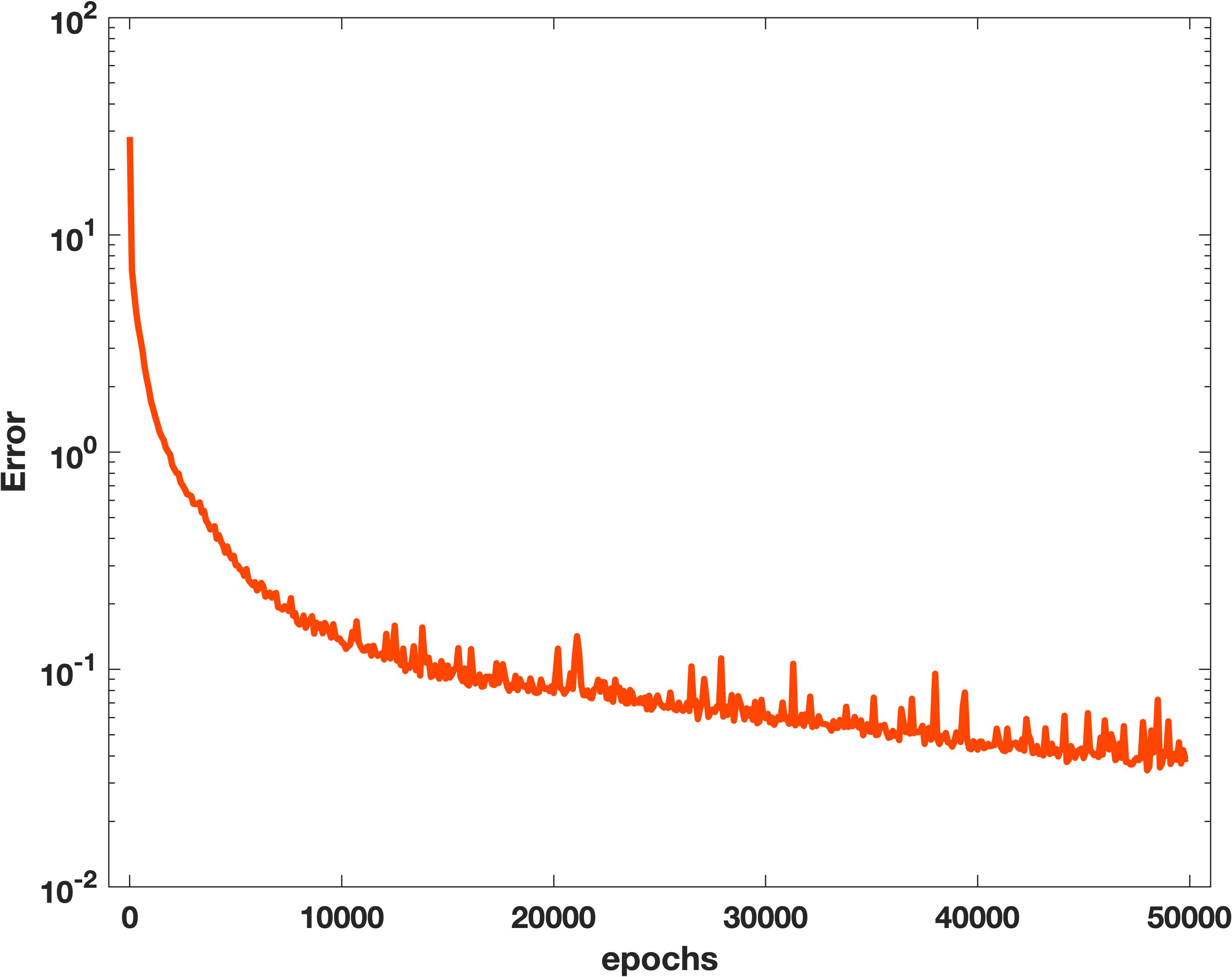}}
   \subcaptionbox{\label{fig:sphere5d_error}}
   {\includegraphics[width=0.47\textwidth]{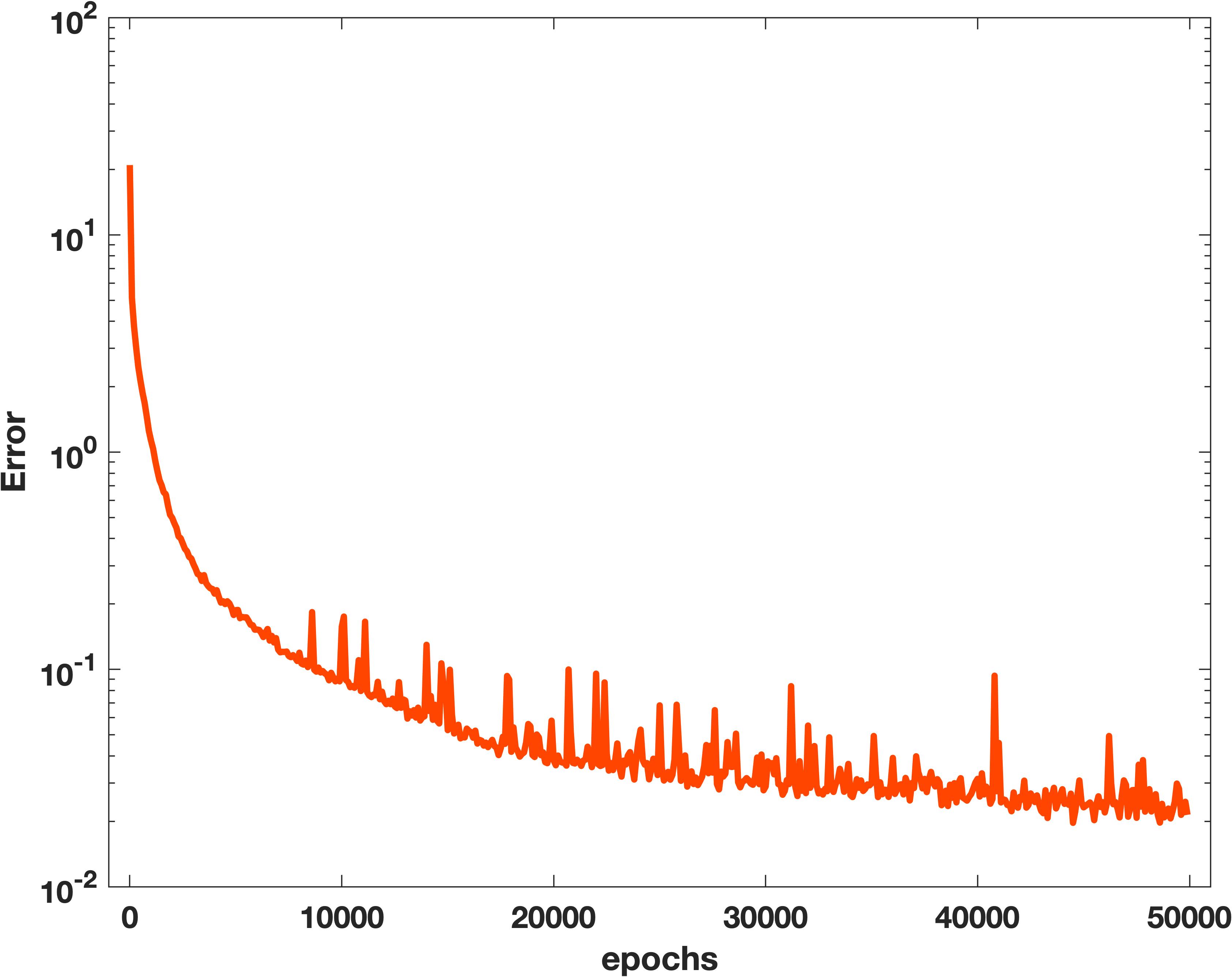}}\\
   \subcaptionbox{\label{fig:sphere10d_error}}
   {\includegraphics[width=0.47\textwidth]{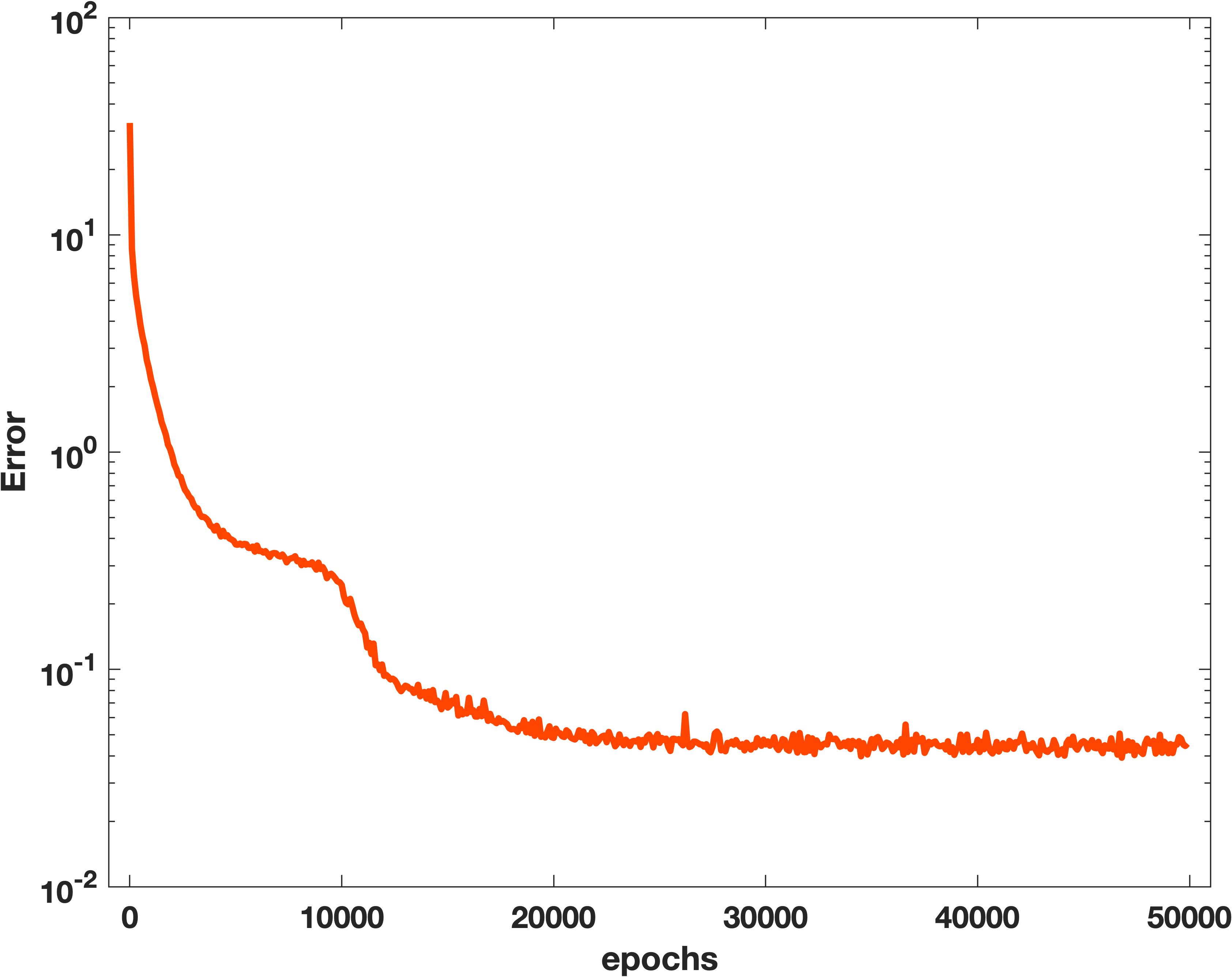}}
   \subcaptionbox{\label{fig:sphere20d_error}}
   {\includegraphics[width=0.47\textwidth]{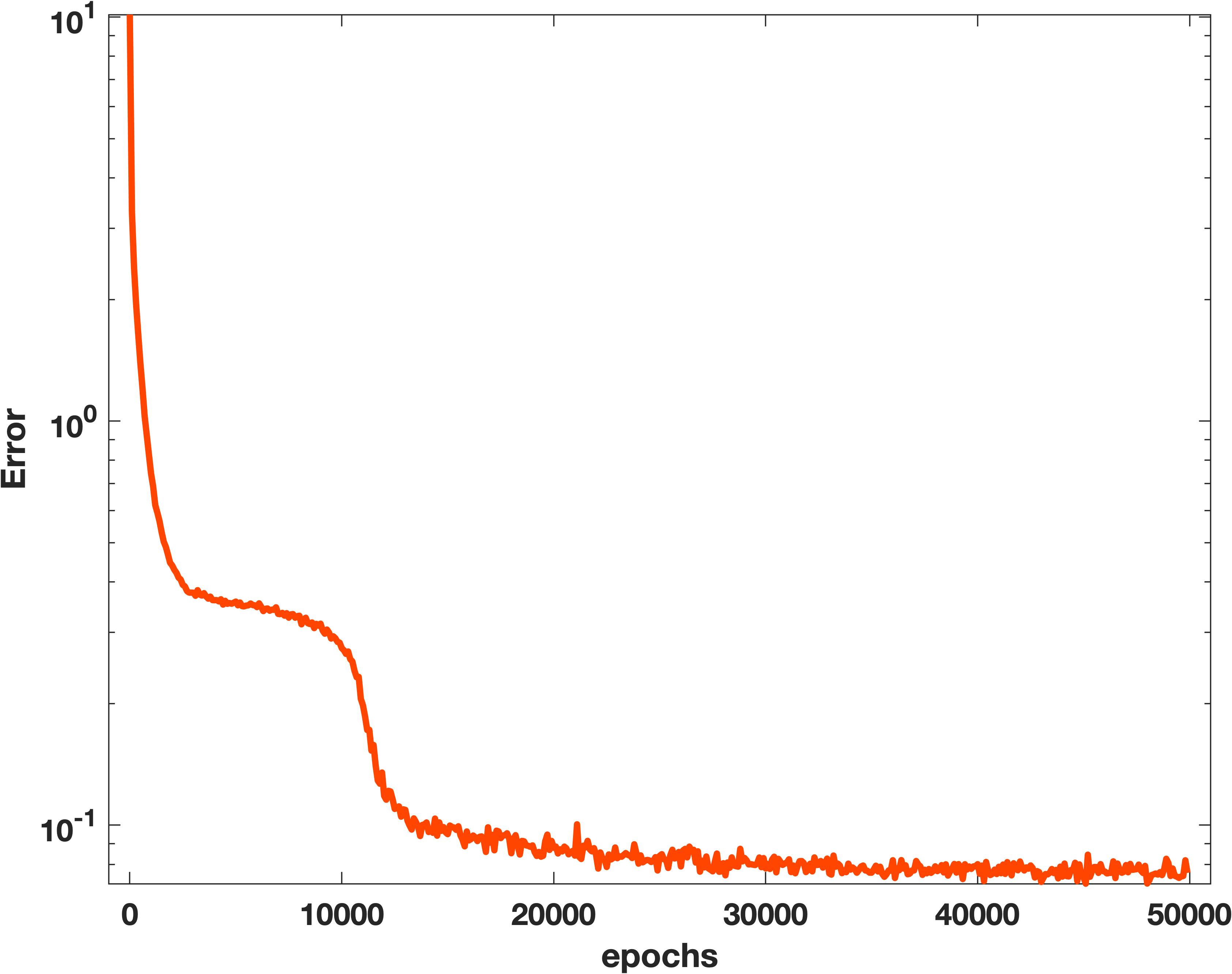}}
   \caption{Decay of loss functions in high-dimensional interface problems: (a) $d=3$; (b) $d=5$; (c) $d=10$; (d) $d=20$. }\label{fig:sphere}
\end{figure}

\begin{table}[h!] 
\centering
 \caption{Relative $L^2$ errors for high dimensional interface problems}\label{tab:highdim}
 \begin{tabular}{|c | c | c | c | c|} 
 \hline
 $d$ & $3$ & $5$ & $10$ & $20$ \\  \hline
Error($\%$) & 3.90 & 2.14 & 4.43 & 7.04  \\
 \hline
 \end{tabular}

\end{table}

In this case, the measure of $\Omega_1 = \frac{\pi^{d/2}r_0^d}{\Gamma(d/2+1)} $ and the measure of $\Omega_2 = 1 - \frac{\pi^{d/2}r_0^d}{\Gamma(d/2+1)}$. The measure of the interface is  $\Gamma=\frac{2\pi^{d/2}r_0^{d-1}}{\Gamma(d/2)}$ . 
We approximate each component of the solution by ResNet with 3 blocks of 20 neurons. In total, there are 2621 parameters for each deep neural network.  To generate uniformly random points on $d$ dimensional spheres, we  use the \verb+sample_hypersphere+ function in BoTorch \cite{BKJD2020}.  To guarantee there are enough random sampled point in $\Omega_1$. We first generate $102$  (about $1/10$ $N=N_1 + N_2$) uniformly random point inside the $d$-dimensional ball with $r_0$ using the dropped coordinates method \cite{HaLa2010}.  In specific, we firstly generate $102$ points on $d+2$ sphere with radius $r_0$ using  the \verb+sample_hypersphere+ function in BoTorch \cite{BKJD2020} and drop the last two coordinates to get the uniform random sampled points. Then, we generate $912$ random points in $\Omega$ which may also contain points in $\Omega_1$. So we have $N_1\ge 102$ and $N_2=1024-N_1$. We take $N_f = 256$ and $N_b = 256d$.  In the following test, we take $\gamma_b = 3000$ and $\gamma_f = 500$.

We consider four different high dimensional cases: $d=3, 5, 10, 20$. Figure \ref{fig:sphere} displays the decay of the relative $L^2$ error during the training process.  Similar to the previous two examples, we can see that the errors decay quickly at the first few epochs and then fluctuate around some specific levels.  In Table \ref{tab:highdim}, we reports the relative $L^2$ error after $50000$ epochs.  We can see that we get almost the same level of relative $L^2$ errors independent of the dimensionality of the space even though we use the same number of points.

 \section{Conclusion}\label{sec:conclusion}
 As a continuous study of our previous works on elliptic interface problems \cite{guo2017gradient,guo2018gradient,GuYa2018}, we propose a deep unfitted Nitsche method to address the high-dimensional challenge with high-contrasts. The challenge is also known as the curse of dimensionality. The classical numerical methods like finite difference and finite methods require unaffordable computational time.  The proposed method deploys the deep neural network to solve the equivalent high-dimensional optimization problem. To address the high contrasts of the solution, we introduced a so-called unfitted Nitsche energy functional, which utilize different deep neural networks to present different components of the solution in the high dimensional case. Different deep networks are patched together weakly by Nitsche's method and can be trained independently using the unfitted Nitsche functional.   The unfitted Nitsche's energy function is approximated by the Monte-Carlo methods. An additional penalty term is added to the discrete energy functional to handle the Dirichlet boundary conditions. The proposed method is easy to be implemented and mesh-free, which is illustrated by several numerical examples including high contrasts and high dimensional cases.


\section*{Acknowledgment}
H.G. is grateful to Prof. Zuoqiang Shi from Tsinghua University and Dr Guozhi Dong from the Humboldt University of Berlin for stimulating and helpful discussions. 
H.G. was partially supported by Andrew Sisson Fund of the University of Melbourne, X.Y. was partially supported by the NSF grants  DMS-1818592 and DMS-2109116.

\bibliographystyle{siamplain}
\bibliography{references}
\end{document}